\documentclass{amsart}
\usepackage{amssymb}
\usepackage{color}
\definecolor{darkgreen}{cmyk}{1,0,1,.2}
\definecolor{m}{rgb}{1,0.1,1}
\definecolor{green}{cmyk}{1,0,1,0}

\definecolor{test}{rgb}{1,0,0}
\definecolor{cmyk}{cmyk}{0,1,1,0}

\newtheorem{Equation}{}[section]

\newtheorem{theorem}[Equation]{Theorem}
\newtheorem{proposition}[Equation]{Proposition}
\newtheorem{lemma}[Equation]{Lemma}
\newtheorem{corollary}[Equation]{Corollary}

\newtheorem{definition}[Equation]{Definition}

\newtheorem{remark}[Equation]{Remark}

\def\ch{\operatorname{ch}}

\def\codim{\operatorname{codim}}

\def\HH{\operatorname{HH}}
\def\ind{\operatorname{ind}}
\def\End{\operatorname{End}}

\def\Iso{\operatorname{Iso}}

\def\Hom{\operatorname{Hom}}

\def\Ind{\operatorname{Ind}}

\def\Td{\operatorname{Td}}

\def\tr{\operatorname{tr}}
\def\sign{\operatorname{sign}}

\def\C{\mathbb C}

\def\R{\mathbb R}

\def\Z{\mathbb Z}

\def\maA{{\mathcal A}}
\def\maB{{\mathcal B}}

\def\maG{{\mathcal G}}
\def\cG{{\mathcal G}}
\def\maH{{\mathcal H}}

\def\maN{{\mathcal N}}
\def\maU{{\mathcal U}}

\def\cP{{\mathcal P}}

\def\what{\widehat}
\def\wtit{\widetilde}

\def\te{\wtit e}

\def\wh{\widetilde{h}}

\def\dd{\displaystyle}
\def\pa{\partial}

\def\ep{\epsilon}

\begin{document}

%%%%%%%%%%%%%%%%%%%%%%%%%%%%%%%%%%%%%%%%%%%%%%%%%%%%%%%
%                                                     %
%   THE MANUSCRIPT BEGINS HERE                        %
%                                                     %
%%%%%%%%%%%%%%%%%%%%%%%%%%%%%%%%%%%%%%%%%%%%%%%%%%%%%%%

%%% TITLE

\title[Higher Lefschetz theorem \today]
{The higher fixed point theorem for foliations\\ I. Holonomy invariant currents\\
\today}

%%% AND THE AUTHORS ARE:

\author[M-T. Benameur]{Moulay-Tahar Benameur}
\address{UMR 7122 du CNRS, Universit\'{e} de Metz, Ile du Saulcy, Metz, France}
\email{benameur@math.univ-metz.fr}
\author[J.  L.  Heitsch \today]{James L.  Heitsch}
\address{Mathematics, Statistics, and Computer Science, University of Illinois at Chicago} 
\email{heitsch@math.uic.edu}
\address{ Mathematics, Northwestern University}

%\email{j-heitsch@northwestern.edu}


\begin{abstract} In this paper, we prove a higher Lefschetz formula for foliations in the presence of a closed Haefliger current. To this end,
we associate with such a current an equivariant cyclic cohomology class of Connes' $C^*$-algebra of the foliation, and compute its pairing with the
localized equivariant $K$-theory in terms of local contributions near the fixed points. 
\end{abstract}

\maketitle
\tableofcontents

\section{Introduction}

Let $V$ be a smooth closed manifold and let $f: V \to V$ be a smooth map. Then  $f$ induces endomorphisms $H^i(f)$ on the finite dimensional $i$-th cohomology spaces $H^i(V, \R)$. The Lefschetz number of $f$ is the integer given by the alternating sum $L(f) := \sum_{i=0}^{\dim (V)} (-1)^i \tr (H^i(f))$, where $\tr$ is the ordinary trace.   It is a topological invariant of $f$.  When $f$ has a finite fixed point set, the classical Lefschetz fixed point formula states
$$
L(f) = \sum_{f(x)=x} \sign \det (I - d_xf).
$$
In \cite{AB}, under appropriate assumptions, Atiyah and Bott extended the above construction to the geometric elliptic complexes over $V$. More precisely, for any elliptic complex $(E,d)$ over $V$, and under appropriate compatibility conditions, Atiyah and Bott defined 
$$
	L(f;E,d)=\sum_{i\geq 0} (-1)^i \tr(f^*|_{H^i(E,d)}) \quad \in \quad  \C,
$$
where $\tr$ is again the ordinary trace of endomorphisms of the finite dimensional $i$-th cohomology space $H^i(E,d)$ of $(E,d)$, and $f^*|_{H^i(E,d)}$ denotes the endomorphism of $H^i(E,d)$ induced by the action of $f$. In general, such a number need not be an integer, or even a real number. When $(E,d)$ is the de Rham complex of differential forms on $V$, the Atiyah-Bott Lefschetz number coincides with the classical topological Lefschetz number $L(f)$ of $f$, using the de Rham isomorphism. Moreover, using heat kernel methods, Atiyah and Bott extended the Lefschetz fixed point formula to such geometric complexes and expressed $L(f;E,d)$ in terms of appropriate topological invariants of the fixed point submanifold. Consequences of such formulae are nowadays well known and encompass results in topology, in complex geometry as well as in number theory, see for instance \cite{HirzebruchZagier}.

The relationship between the Atiyah-Bott Lefschetz formulae and index theory should be transparent from the very definition of $L(f;E,d)$. More precisely, assume that $f$ is an isometry for some metric $g$ on $V$ and denote by $H$ the compact Lie group generated  by $f$ in  the group $\Iso(V,g)$ of isometries of $(V,g)$. Then for any $H$-invariant elliptic \textit{pseudodifferential} complex $(E,d)$, the Lefschetz number $L(f;E,d)$ is well defined and can be re-interpreted as the evaluation at $f$ of a character, namely the $H$-equivariant analytical index of the complex $(E,d)$. In \cite{AS2}, Atiyah and Segal gave a proof of the Lefschetz fixed point formulae with respect to pseudodifferential complexes, as a consequence of the $H$-equivariant index theorem of Atiyah-Singer. Moreover, the index method allowed them  to extend the fixed point formulae to higher dimensional fixed point sets with interesting new characteristic invariants \cite{AS1}.\\

When the closed manifold $V$ is endowed with a smooth foliation $F$, the Lefschetz problem can be stated for leaf-preserving maps and for leafwise elliptic pseudodifferential complexes. The Atiyah-Bott heat kernel approach for the geometric complexes has  been extended in \cite{HL1} to foliations which admit a holonomy invariant transverse measure $\Lambda$. In this case, the Lefschetz number of a leaf-preserving diffeomorphism $f$ with respect to the leafwise elliptic complex $(E,d)$ is defined by replacing the usual trace $\tr$ by a von Neumann trace $\tr_\Lambda$ associated with $\Lambda$, and the precise definition is as expected
$$
L_\Lambda (f; E,d) := \sum_{i\geq 0} (-1)^i \tr_\Lambda (f^*|_{H^i(E,d)}) \quad \in \quad  \C.
$$
Moreover, under the usual transversality assumption, the measured MacKean-Singer formula allowed the proof of a measured Lefschetz fixed point theorem, see again \cite{HL1}. When $f$ is an isometry of $V$ for some metric $g$, it was proved in \cite{Be1} that the measured Lefschetz number can again be interpreted as the evaluation at $f$ of an appropriate character, namely the equivariant measured analytic index introduced by Connes \cite{CInt}. Moreover, the measured Lefschetz fixed point formula can be deduced from an equivariant measured index theorem for foliations exactly as in the Atiyah-Segal approach, see again \cite{Be1}.

One of the main features of Connes' approach to foliation index theory is that the measured index of a leafwise elliptic complex is the von Neumann trace of a much more fundamental object: the $K$-theory index in the $C^*$-algebra of the foliation \cite{CInt}. This latter exists even if the foliation has no holonomy invariant measure and yields more sophisticated index invariants by the use of Connes' cyclic theory \cite{C85}. In this picture, holonomy invariant measures correspond to the pairing of the $K$-theory index class with a zero degree cyclic cocycle.  Moreover, one can use the transverse geometry of the foliation to produce higher degree cyclic cocycles and higher indices. In \cite{C86} for instance, examples of cyclic cocycles arising from the so-called Haefliger currents are used to investigate higher index theory and some of its deep consequences.  On the other hand, the main result of \cite{Be1} is a proof of a $K$-theory Lefschetz fixed point formula which takes place in the, localized with respect to $f$, equivariant $K$-theory of Connes' $C^*$-algebra of the foliation. As recalled in Section \ref{Ktheory}, this formula is a consequence of the $K$-theory index theorem proved in \cite{CS}.  In order to extend the measured Lefschetz theorem to higher degree cyclic cocycles, one is naturally led to the study of the following interesting questions for a given compact Lie group $H$ of leaf-preserving (or more generally $TF$-preserving) diffeomorphisms of $V$:
\begin{enumerate}
 \item Show that the cyclic cocycles arising from the Haefliger homology of the foliation are $H$-equivariant cyclic cocycles.
\item Show that evaluation at a conjugacy class $[h]$ of $H$ of the \underline{well defined} pairing of these $H$-equivariant cocycles with the equivariant $K$-theory of Connes' $C^*$-algebra induces a well defined pairing with the localized $K$-theory at the prime ideal in $R(H)$ associated with $[h]$.
\item Once $(1)$ and $(2)$ have been accomplished, show that  the equivariant pairing of Haefliger homology with $K$-theory is compatible with shriek maps associated with transverse submanifolds and their foliations.
\end{enumerate}
We give here a complete solution to the three problems. The first one allows us to define the higher Lefschetz number of any leaf-preserving map $f:V\to V$ which is isometric for some metric on $V$. The second and third ones allow us to prove our higher Lefschetz theorem. It is worth pointing out that the answer to the third problem allows us to prove a higher Lefschetz formula which  only depends on the fixed point submanifold with its {\em induced} foliation, which is in the spirit of the classical Lefschetz theorems. Indeed, a formula depending only on the saturation of this fixed point submanifold can be obtained in a much easier way and for a larger class of cyclic cocycles.   However,  it doesn't really qualify as a fixed point formula, as it would not be useful in most of the examples where this saturation coincides with the ambiant manifold. 
We also translate in this paper the higher fixed point formula to cohomology by using the Chern-Connes character and the Grothendieck-Riemann-Roch theorem proved in \cite{BHI}. 

To sum up,  we use geometric \textit{equivariant} cyclic cocycles of arbitrary even dimension to produce Lefschetz fixed point formulae in terms of characteristic classes at the fixed point submanifold with its induced foliation. When the foliation is top dimensional, we recover the results of Atiyah-Segal \cite{AS2} and Atiyah-Singer \cite{AS1}. When the foliation admits a transverse holonomy invariant measure, we recover the results of \cite{HL1}. When the foliation is given by the connected components of the fibers of a smooth closed fibration, we recover the results obtained by the first author in \cite{Be3}. Finally, since we may use any closed Haefliger current on the foliation, we  get new and very interesting formulae.

We now explain more explicitely our results. Denote by $h: V,F \to V,F$ a leaf-preserving map, and assume that it generates a compact Lie group $H$ in some Lie group of maps of $V$. 
So, we do not assume a priori that $h$ a leafwise isometry.  Note also that $H$ need not preserve the leaves although it will always send leaves to leaves.  The Chern-Connes construction allows us to construct out of any closed Haefliger current $C$, an equivariant
cyclic cocycle $\tau_C$ which pairs with the equivariant $K$-theory of the foliation and gives a central function on the compact Lie group $H$. This is the content of 

\medskip\noindent
{\bf{Theorem}} \ref{locp}
{\em{  Let $C$ be an even dimensional closed Haefliger current on $(V, F)$ and let $\tau_C$ be the associated cyclic cocycle. Then for any $H$-vector bundle $E$, $\tau_C$ is an equivariant cyclic cocycle on $C_c^{\infty}(\cG, E)$. Moreover, the equivariant pairing extends  to a well defined pairing
 $$
K^H (C^*(V,F))\otimes H_{ev}(V/F) \longrightarrow C(H)^H.
$$
}}
The pairing of $\tau_C$ with $\Ind^H (E,d)$, (the equivariant index class  in $K$-theory of a leafwise elliptic pseudodifferential complex $(E,d)$), evaluated at $h$ will be, by definition, our higher Lefschetz number $L_C (h; E,d)$:
$$
L_C (h; E,d) := <\tau_C, \Ind^H (E,d) > (h) \quad \in \quad \C.
$$

Compatibility with the prime ideal associated with $h$ is then easily deduced from the properties of the equivariant pairing. Finally, the commutation with transverse shriek maps is  proved in Proposition \ref{trans}. The higher Lefschetz formula that we get using the higher index theorem of \cite{BHI} is then stated as follows when the fixed point submanifold $V^h$ is transverse to the foliation with an oriented induced foliation $F^h$:

\medskip\noindent
{\bf{Theorem}} \ref{basic2}  
{\em{For any even dimensional closed Haefliger current $C$, 
$$
	L_C(h; E,d)=\left< \int_{F^h} \frac{\ch_{\C}(i^{*}[\sigma(E,d)](h))}{\ch_{\C}(\lambda_{-1}(N^{h}
	\otimes{\C})(h))} \, \Td(TF^{h} \otimes \C) , C|_{V^h} \right >.
$$}}

\noindent
See Section \ref{HigherLef} for more details and for the notation.
Theorem \ref{basic2} simplifies notably when the fixed point 
submanifold $V^{h}$ is a strict transversal, which corresponds in the case of a foliation by a single leaf  to 
the original case of isolated fixed points:

\medskip\noindent
{\bf{Theorem}} \ref{sttrans}
{\em{If $V^h$ is a strict transversal, }}
$$
L_C(h;E,d)  \,\, = \,\,  \left < \frac{ \sum_{i} (-1)^{i} \ch_{\C}([E^{i}|_{V^{h}}](h))}{\sum_{j} (-1)^{j} \ch_{\C}
([\wedge^{j}(TF|_{V^{h}} \otimes \C)](h))} , C|_{V^h} \right >.
$$

\noindent
As an obvious corollary of Theorem \ref{basic2}, we see for instance that the non vanishing of the higher Lefschetz numbers, for some current $C$, immediately implies the existence of fixed points under the action of $h$. 

In order to keep the present paper to a reasonable size, we have postponed the investigation of the applications of our higher Lefschetz theorem to complex geometry, topology and number theory to the second part \cite{BHV}. In particular, we deduce in \cite{BHV} a higher version of the Atiyah-Hirzebruch rigidity theorems, new formulae for appropriate sums of products of cotangents, as well as some integrality consequences for Riemannian foliations.

We now describe the contents of each section. Section 2 gives a brief survey of the $K$-theory Lefschetz theorem. In section 3, we define the appropriate equivariant cyclic cohomology and construct a pairing between this equivariant cyclic cohomology and equivariant $K$-theory which factors through the representation ring of the compact Lie group generated by the isometry $f$. In section 4, we return to the foliation case and construct a  correspondence between Haefliger homology and equivariant cyclic cohomology.  Finally, in section 5 we use the results of the previous sections together with the $K$-theory Lefschetz theorem to prove the higher Lefschetz formula.\\

\noindent
{\em Acknowledgements.} The authors wish to thank  Alain Connes, Thierry Fack, Gilbert Hector, Steve Hurder, Victor Nistor, and Denis Perrot for  helpful discussions. They are also indebted to the referee for his judicious comments.

\section{Review of the $K$-theory Lefschetz theorem.}\label{Ktheory}

Most of the results of this paper are based on the $K$-theory Lefschetz theorem proven by the first author in \cite{Be1} and which we now  recall.   Let $F$ be a smooth foliation of the smooth compact Riemannian manifold $(V,g)$ and denote the tangent bundle to the foliation $F$ by $TF$, and its normal bundle by $\nu$. Denote by $\cG$ the holonomy groupoid of $F$, which consists of equivalence classes of leafwise paths, where two paths are identified if they start at the same point, end at the same point, and the holonomy germ along them is the same.  Composition of paths makes $\cG$ a groupoid, and its space of units $\cG^0$ consists of the classes of the constant paths, so  $\cG^0 \simeq V$.  Denote  
by $\cG_x$ the elements of $\cG$ which start at the point $x \in V$,  by $\cG^y$ those elements  which end at the point $y \in V$, and  by $\cG^y_x$ the intersection $ \cG_x \cap \cG^y$.  
We have the maps $s,r:\cG \to V$, where $s(\gamma) = x$, if $\gamma  \in \cG_x$, and $r(\gamma) = y$ of $\gamma \in \cG^y$.  The metric $g$ on $V$ induces a canonical metric on $\cG$, and so the splitting $T\cG = TF_s \oplus TF_r \oplus \nu_{\cG}$.  Note that  $r_*(\nu_{\cG,\gamma}) = \nu_{r (\gamma)}$, and $s_*(\nu_{\cG,\gamma}) = \nu_{s (\gamma)}$.
For details, see \cite {BHI}. 
The metric on $\cG$ gives metrics on the  submanifolds $\cG_x, \cG^y \subset \cG$.   So objects such as $L^{2}(\cG_{x})$ and $L^{2}(\cG^y)$ are well defined, and do not depend on the choice of metric since $V$ is compact.  
Note that $r:\cG_{x} \to L_x$ is the holonomy covering of $L_x$, the leaf of $F$ through $x$, and similarly, $s:\cG^y \to L_y$ is  the holonomy covering of $L_y$.

Let $H \subset \Iso (V,g)$ be a compact Lie group which acts by  $TF$-preserving isometries on $(V,g)$.  The isometry of $(V,g)$ which corresponds to the action of $h\in H$ will be denoted by $h$ for simplicity, so such $h$ takes leaves of $F$ to leaves of $F$, but does not necessarily take leaves to themselves. An easy example is the action of the torus on its constant slope foliation. When the slope is irrational, this action can be seen as being topologically generated by a leaf-preserving isometry, but the whole group does not preserve the leaves.  
Connes' $C^*$-algebra of $(V,F)$  is denoted $C^*(V,F)$, see  \cite{C82}.  It is easy to check that $C^*(V,F)$ is an $H$-algebra, i.e. the induced action of $H$ is strongly continuous for the $C^*$-norm.
Let $P$ be an $H$-invariant, uniformly supported, elliptic  pseudodifferential $\cG$-operator acting from sections of the $H$-vector bundle $E^+$ to sections of the $H$-vector bundle $E^-$. See  \cite{CInt} for the precise classical definitions.   Choose an $H$-invariant Hermitian structure on $E$ and denote by $\epsilon_{V,E}$ the Hilbert $C^*$-module over $C^*(V,F)$, associated with the continous field of Hilbert spaces $L^{2}(\cG_{x}, r^* E)$.   Then $(\varepsilon_{V,E}, P)$ defines a $KK$-class  in the Kasparov equivariant
group $KK^{H}(\C, C^*(V,F))$,  \cite{Kasparov80}. The image of this class under the isomorphism 
$$
KK^{H}(\C, C^*(V,F)) \simeq K^H( C^*(V,F)),
$$
is  the analytic  $H$-index of $P$, and is denoted $\Ind^{H}_{a,V} (P)$. So, $\Ind^{H}_{a,V} (P) \in K^H( C^*(V,F))$.

Denote the space of leaves of $F$ by $V/F$, and recall that a map $g:N \to V/F$, from a smooth manifold $N$ to the space of leaves, is by definition a $\cG$-valued 1-cocycle over $N$, see \cite{CS}.   
Recall the $H$-submersion 
$$
	p : TF \rightarrow V/F
$$ 
which is the composite map $TF \to V \to V/F$. 

The pull back by $p$ of the  ``tangent bundle"  $T(V/F)$ is well defined, and it is just the pull-back $\pi^*\nu$ of the normal bundle $\nu = TF^{\perp} \subset TV$ of the foliation $F$ by the projection $\pi: TF \to V$.   The map $p$ is  $K^H$-oriented, which means that the vector bundle $T(TF) \oplus \pi^*\nu$ admits a spin$^c$ structure, which is given by an $H$-equivariant Hermitian vector bundle of irreducible representations of the Clifford algebra bundle associated  with $TF \oplus \pi^*\nu$.  
Using this spin$^c$-bundle, Connes and Skandalis constructed in \cite{CS} a Gysin class  
$$
	p_! \in KK_{H}(C_0(TF),C^*(V,F)),
$$
which generalizes the Atiyah-Singer topological shrieck maps \cite{AS1}.
The topological $H$-index of $P$ is by definition
$$
	\Ind^{H}_{t,V} (P) = [\sigma(P)]\otimes_{TF} p_! \quad \in \;\; K^H (C^*(V,F)),
$$ 
where $[\sigma(P)]$ is the
$H$-equivariant class of the
principal symbol of $P$ in $K_{H}(TF)$, and we have used the Kasparov product over the $C^*$-algebra $C_0(TF)$ together with the isomorphism $KK_{H}(\C,C^*(V,F)) \simeq K^H (C^*(V,F))$ to see the $H$-index as an equivariant $K$-theory class. 
The leafwise $H$-equivariant index theorem of \cite{C85} is then 

\begin{theorem}\   
$$
\Ind^{H}_{a,V}(P) = \Ind^{H}_{t,V}(P)\quad \in \;\; K^H(C^*(V,F)).
$$
\end{theorem}

This theorem clearly works as well for elliptic $H$-invariant $\cG$-complexes, see \cite{AS1}.\\

We  now to state the $K$-theory Lefschetz theorem for a fixed isometry $f\in H$.
Let $(E,d)$ be an elliptic pseudodifferential $\cG$-complex   over $(V,F)$.   Recall that such a complex is a finite collection of smooth vector bundles $(E^i)_{0\leq i \leq k}$ over the ambiant manifold $V$ together with a collection of $\cG$-invariant \underline{uniformly supported} pseudodifferential operators $d^i = (d^i_x)_{x\in V}$ acting from the smooth sections of $r^*E^i$ to the smooth sections of $r^*E^{i+1}$, see \cite{CInt, NistorWeinsteinXu}. So, we have
$$
0 \to C_c^\infty (\cG_x, r^*E^0) \stackrel{d^0_x}{\longrightarrow } C_c^\infty (\cG_x, r^*E^1)\stackrel{d^1_x}{\longrightarrow } \cdots \stackrel{d^{k-1}_x}{\longrightarrow }C_c^\infty (\cG_x, r^*E^k) \to 0,
$$
with $d^{i+1}\circ d^i =0$. Notice that the uniform support of $d^i$ means that its Schwartz kernel, viewed as a section over $\cG$ using the $\cG$-invariance, is compactly supported.  Associated to $(E,d)$ is the vector bundle complex of principal symbols that we denote by $\sigma(E,d)$, see again \cite{CInt},
$$
0 \to \pi^*E^0 \stackrel{\sigma(d^0)}{\longrightarrow } \pi^*E^1 \stackrel{\sigma(d^1)}{\longrightarrow } \cdots \stackrel{\sigma(d^{k-1})}{\longrightarrow }\pi^*E^k \to 0,
$$
where $\pi: TF\to V$ is the projection. Ellipticity of $(E,d)$ then means that this principal symbol complex is exact off the zero section of $(TF, \pi, V)$. We assume that $(E,d)$  is  $H$-invariant in the sense that each $E^i$ is an $H$-equivariant vector bundle over the $H$-manifold $V$, and the operators $d^i$ commute with the naturally induced actions of $H$, \cite{Be1}.

Let $I_{[f]} = \{\chi \in R(H) \, | \, \chi (gfg^{-1}) = 0, \forall g\in H\}$ be the prime ideal,  in the representation ring $R(H)$ of $H$, associated with the conjugacy class $[f]$ of $f$. Localization of the ring $R(H)$ with respect to $I_{[f]}$ yields the ring of fractions that we denote as usual by  $R(H)_{[f]}$. Given an $R(H)$-module $M$, we denote by $M_{[f]}$ the localization of $M$ 
 with respect to $I_{[f]}$.   So, $M_{[f]}$ is an $R(H)_{[f]}$-module. When $H$ is abelian, we remove the brackets and denote as well by $f$ the conjugacy class of $f$. Important examples of such modules are the $H$-equivariant $K$-theories of $H$-algebras with the module structure given by Kasparov product. We now  define the Lefschetz class as in  \cite{Be1}.

\begin{definition}\label{fclass} \begin{enumerate}
 \item The Lefschetz class $L(f;E,d)$ of the conjugacy class $[f]$ of the isometry $f$, with respect to the $H$-invariant elliptic $\cG$-complex $(E,d)$ is the localized analytic $H$-index of $(E,d)$ with respect to the ideal $I_{[f]}$.   More precisely,
$$
	L([f] ; E,d) \,\, = \,\, \frac{\Ind^{H}_{a,V}(E,d)}{1_{R(H)}} \quad \in \;\; K^{H}(C^*(V,F))_{[f]}.
$$
\item Denote by $H_1$ the compact  subgroup of $H$ generated by $f$.   The Lefschetz class of $f$ with respect to the $H_1$-invariant elliptic $\cG$-complex $(E,d)$ is the localized analytic $H_1$-index of $(E,d)$ with respect to the ideal $I_{f}$ of $R(H_1)$.   More precisely,
$$
	L(f ; E,d) \,\, = \,\, \frac{\Ind^{H_1}_{a,V}(E,d)}{1_{R(H_1)}}  \quad \in \;\;  K^{H_1}(C^*(V,F))_{f}.
$$
                                 \end{enumerate}

\end{definition}

One can express $L([f] ; E,d)$ in terms of the union of the fixed points of all the elements of the conjugacy class $[f]$ as in \cite{AS2}. Notice now that the usual Lefschetz formulae express Lefschetz numbers of a given map in terms of local contributions from the fixed point submanifold of that map. In our context of the foliated Lefschetz problem, this means that we can concentrate on the compact group $H_1$ topologically generated by a given $h\in H$. We shall therefore assume for the rest of this section that $H_1=H$ and only consider $L(f;E,d)$. 

For the foliation with one leaf, namely the manifold $V$ itself, $(E,d)$ becomes
a classical pseudodifferential elliptic $H$-invariant complex over the compact manifold $V$, and
the $R(H)_f$-module $K^{H}(C^*(V,F))_{f}$ coincides with $K^{H}({\mathcal K}(L^2(V))_{f}$,
where ${\mathcal K}(L^2(V))$ is the elementary $C^*$-algebra of compact operators on the Hilbert space $L^2(V)$ of $L^2$ functions on $V$.  Given any smooth kernel $k\in C^\infty (V\times V)$, we define for any $h\in H$ the smooth kernel $k^h(x,y) =k(h^{-1}x,y)$.  The integral of $k^h$ over the diagonal of $V\times V$ then realizes  the Morita equivalence which induces the isomorphism $K^{H}({\mathcal K}(L^2(V))\cong R(H)$.  
 The fixed point submanifold of $f$ is denoted $V^{f}$. We are only interested in the case where $V^{f}$ is transverse to the foliation (this is true when $H$ is connected and preserves the leaves, see \cite{HL1}). Let  $F^{f}$ be the foliation of $V^{f}$ whose tangent bundle is $TF^{f} = TV^{f}\cap TF$, where   $(TV)^{f} = T(V^{f})$ is the fixed point tangent manifold. 
Let $i : TF^{f}\hookrightarrow TF$ be the inclusion. Then an easy inspection shows that $i$ is a $K^H$-oriented map. Thus, following \cite{AS1}, one defines a shrieck element  $i_! \in KK_{H}(TF^{f},TF)$,  \cite{Be1}. 

Notice that since the action of $H$ on $(V^f, F^f)$ is trivial, the $H$-equivariant 
$K$-theory of the manifold $F^f$ can be shown to be  isomorphic to the tensor product $R(H)$-module $K(F^f) \otimes R(H)$. The same isomorphism holds for the $C^*$-algebra of the foliation $(V^f, F^f)$.  Moreover,  the $H$-equivariant index map for the foliated manifold $(V^f,F^f)$ then respects this isomorphism, i.e.
$$
	Ind_{V^f}^H \cong Ind_{V^f} \otimes R(H) : K_H(F^f) \cong K(F^f)\otimes R(H) \to 
	K^H(C^*(V^f,F^f))\cong K(C^*(V^f,F^f))\otimes R(H).
$$

\begin{theorem}\label{abgrp} \cite{Be1} With the above notations, the following diagram commutes 

\hspace{0.25in}
\begin{picture}(415,80)

\put(75,60){$ K(F^{f}) \otimes R(H)\cong K_{H}(F^{f})$}
\put(168,50){ $\vector(4,-3){45}$}
\put(75,35){$Ind_{V^{f}} \otimes R(H) \cong \Ind^{H}_{V^f}$}

\put(215,70){$i_!$}
\put(200,64){\vector(1,0){35}}

\put(243,60){$K_{H}(F)$}
\put(255,50){ $\vector(0,-1){30}$}
\put(125,5){$K(C^*(V^f,F^f))\otimes R(H) \cong K^{H}(C^*(V,F)).$}
\put(265,35){$ \Ind^{H}_{V}$}

\end{picture}
\end{theorem}

Here we have denoted by $\Ind^{H}_{V^f}$ the analytic ($=$ topological) leafwise $H$-index for the compact foliated manifold $(V^f,F^f)$. This index map takes values in the $R(H)$-module $K^H(C^*(V^f,F^f))$ and we have implicitely used composition with a quasi trivial element of $KK^H(C^*(V^f,F^f),C^*(V,F))$ corresponding to  Morita extension as in \cite{C85}, to view its range in $K^{H}(C^*(V,F))$. Now,  functoriality of the localization at the prime ideal $I_f$ yields the $K$-theory Lefschetz Theorem.

\begin{theorem}\label{Lefschetz}\cite{Be1} [The $K$-theory Lefschetz Theorem]
Let $(V,F)$ be a compact foliated manifold and let $f$ be an isometry of $V$ for some metric $g$. Assume that $f$ preserves the  leafwise bundle  $TF$. Then,  the Lefschetz class $L(f;E,d) \in K^H(C^*(V,F))_f$
of $f$ with respect to an $f$-invariant elliptic $\cG$-complex $(E,d)$ is given by
$$
	L(f ; E,d) \,\, = \,\,
	(Ind_{V^{f}} \otimes R(H)_{f}) \left(\frac{i^{\ast}[\sigma(E,d)]}{\lambda_{-1} (N^{f} \otimes \C)}\right),
$$
where $i^*:K_H(F)_f \to K_H(F^f)_f$ is the restriction homomorphism, $N^{f}$ is
the normal $H$-vector bundle to $V^{f}$ in $V$, and $\lambda_{-1}
(N^{f}\otimes \C) = \sum(-1)^{i}
[\wedge^{i}(N^{f} \otimes \C)]$ in $K_{H}(V^{f})_{f} \cong K(V^{f}) \otimes R(H)_{f}.$
\end{theorem}

The reader should note that we do not assume that $f$ preserves the leaves.
That the fraction in the previous theorem is well defined is clear since 
$$
\lambda_{-1}(N^{f}\otimes \C)\,\, := \,\, \sum_i (-1)^i [\wedge^i (N^f\otimes \C)] \quad  \in \;\; K_H (V^f),
$$ 
becomes a unit when pushed to the localized ring $K_{H}(V^{f})_{f}$ and we have used the $K_{H}(V^{f})_{f}$-module structure of $K_H(F^f)_f$. 

Hence, this theorem tells us that the Lefschetz class
$L(f ; E,d)$ coincides with the  index  of a virtual leafwise operator on the fixed point foliation.
When $V^{f}$ is a strict transversal, i.\ e.\ when $V^f$ is transverse to the foliation with dimension equal to
the codimension of the foliation,  we get
$$
	L(f;E,d) = \frac{\sum(-1)^{i}[E^{i}|_{V^f}]}{\sum(-1)^j[\wedge^{j}(F|_{V^f} \otimes \C)]}
	\otimes_{V^{f}}[[V^{f}]]\quad \in \;\;  K^H(C^*(V,F))_f.
$$
Here $[[V^{f}]]$ is the class associated to the transversal $V^{f}$ in $KK^H(C(V^f),$ $C^*(V,F))$ 
as defined in \cite{C82}.  Since Kasparov product by $[[V^{f}]]$ is {\it{$K$-integration over $V^{f}$}},
 the above formula agrees with the classical ones. In particular when $(E,d)$ is  the de Rham complex along the leaves, our Lefschetz class coincides with the {\it{$K$-volume}} of $V^{f}$, a class $[V^{f}]$ represented by an $H$-invariant idempotent in $C_c^{\infty}(\cG)$, supported in a small tubular neighborhood of $V^f$, \cite{C82}.

It is worth pointing out that the above $K$-theory Lefschetz theorem implies the measured Lefschetz theorem of \cite{HL1} in the isometric case. More precisely, if the foliation $(V,F)$ admits a holonomy invariant transverse measure $\Lambda$ \cite{Plante}, then recall the  trace $\tau^\Lambda$ constructed by Connes on the von Neumann algebra of the foliation \cite{CInt}. In \cite{Be1}, the first author constructed an additive map $\tau_\Lambda^f: K^H (C^*(V,F))_f \to \C$  associated with $\Lambda$, such that
$$
\tau_\Lambda^f (L(f;E,d))  \,\, = \,\, L_\Lambda (f;E,d),
$$ 
where $L_\Lambda (f;E,d)$ is the measured Lefschetz number defined in \cite{HL1}. Hence composing Theorem \ref{Lefschetz} with $\tau_\Lambda^f$, the Heitsch-Lazarov measured Lefschetz formula is deduced.

We end this review section by  pointing out that the analytic $H$-index $Ind_V^H(E,d)$ is the image of an element $\ind_V^H (E,d)$ of the equivariant $K$-theory group $K^H(C_c^{\infty}(\cG, E))$ of the algebra of compactly supported smooth sections of the bundle $\Hom (E)$ over the graph $\cG$, whose fiber at $\gamma\in \cG$ is $\Hom(E_{s(\gamma)}, E_{r(\gamma)})$.

\section{Equivariant cyclic cohomology}\label{Equiv}

In this section we recall  the definitions and properties of the equivariant cohomologies for actions of compact groups. References for these notions go back to  \cite{BrylinskiPreprint} where equivariant homology arises as the $E^2$ term of a spectral sequence that computes the homology of a topological crossed product. In the papers \cite{BlockGetzler} and \cite{Gong}, the equivariant homologies are computed for algebras of functions on a smooth manifold, extending to the equivariant setting the HKR theorem of Connes \cite{C85}. The computation of these equivariant (co)homologies for more general algebras with appropriate topologies (more precisely bornologies) and endowed with strongly continuous actions of compact groups was recently achieved in the PhD thesis of C. Voigt \cite{Voigt}. We also point out the earlier computation in \cite{NistorInv}.

Let $H$ be a  compact group which acts continuously on the locally convex (unital) algebra $\maA$. We proceed to define the equivariant cohomology of the dynamical system $(\maA, H)$. Denote by $C (H)$ the continuous $\C$ valued functions on $H$.
The equivariant Hochschild complex $(C^*(\maA, H), b)$  is defined as follows.
The cochains $C^n(\maA, H)$ consist of the continuous functions
$f: \maA^{\otimes_{n+1}} \to C (H)$ such that 
$$
f(h a^0, \cdots, ha^n) (h g h^{-1}) = f(a^0, \cdots, a^n)(g), \quad  \forall g,h\in H, \forall a^j\in \maA .
$$
We denote by $f(a^0, \cdots, a^n|h)$ the scalar $f(a^0, \cdots, a^n)(h)$ for $f\in C^n(\maA, H)$. The equivariant Hochschild differential $b: C^n(\maA, H) \to C^{n+1} (\maA, H)$ is defined by
$$
(bf) (a^0, \cdots, a^{n+1} | h):= (b' f)(a^0, \cdots, a^{n+1} | h) + (-1)^{n+1} f ( h^{-1} (a^{n+1}) a^0  , a^1, \cdots, a^n | h) 
$$
where 
$$
(b' f)(a^0, \cdots, a^{n+1} | h) := \sum_{j=0}^n (-1)^j f( a^0, \cdots, a^j a^{j+1}, \cdots, a^{n+1}| h).
$$
Then one checks that $b'^2=0$ and $b^2=0$.

\begin{definition}
The equivariant Hochschild cohomology of the pair $(\maA, H)$ is  the homology $\HH^*(\maA, H)$ of the equivariant Hochschild complex $(C^*(\maA, H), b)$.
\end{definition}

\begin{definition}
An equivariant Hochschild cochain $f$ is cyclic if the following relation holds
$$
f(a^n, a^0, \cdots, a^{n-1} |h) = (-1)^n f(a^0, \cdots, a^{n-1}, h(a^n) |h). 
$$
\end{definition}
We denote by $
\lambda_H (a^0, \cdots, a^n |h) := (h^{-1} (a^n), a^0, \cdots, a^{n-1}  |h)$
the equivariant  permutation. So, an $n$-cochain $f$ is cyclic if $\lambda_H^*f = (-1)^n f$. We denote by $A$ the equivariant cyclic permutation of cochains
$$
A (f)  = \sum_{j=0}^{n-1} (-1)^{nj} (\lambda_H^*)^j (f).
$$
The subspace of $C^*(\maA, H)$ composed of equivariant cyclic cochains is denoted $C^*_\lambda (\maA, H)$. An easy computation shows that 
$$
A\circ b' = b\circ A.
$$
So, the subspace $C^*_\lambda (\maA, H)$ is preserved by $b$ and the  subcomplex $(C^*_\lambda (\maA, H), b)$ of the equivariant Hochschild complex is called the equivariant cyclic complex. Its homology $H_\lambda^* (\maA, H)$ is the equivariant cyclic cohomology of the pair $(\maA, H)$.

\begin{remark}
When $\maA$ is complete, a chain equivalence between the equivariant Hochschild (resp. cyclic) complex of $(\maA, H)$ and the Hochschild (resp. cyclic) complex of the crossed product algebra $\maA\rtimes H$, is constructed in \cite{BT}.   When the group $H$ is in addition finite, this chain map  turns out to be an isomorphism between the equivariant cohomologies and the {\em twisted} cohomologies of the crossed product, as constructed in \cite{BenameurFoliated}. These latter complexes being other alternatives for describing the Hochschild and cyclic cohomologies of the crossed product. 
\end{remark}

For the algebra of smooth functions on a closed manifold, equivariant Hochschild (resp. cyclic, resp. periodic cyclic) (co)homology has been computed in \cite{BlockGetzler} using Borel (co)homology.

Now let $H \to \maU(X)$ be a given finite dimensional unitary representation of $H$ in the Hermitian vector space $X$.  Then $H$ acts on $\maA \otimes \End(X)$ according to the formula
$$
h (a\otimes A) := ha \otimes U(h)AU(h^{-1}), \quad a\in \maA, A\in \End(X) \text{ and } h\in H.
$$

If $f$ is a continuous  equivariant $k$-cochain on the pair $(\maA, H)$, then we define the cochain $\tau\sharp \tr$  on $(\maA\otimes \End(X), H)$ by setting
$$
(f \sharp \tr) (a^0\otimes A^0, \cdots, a^k \otimes A^k | h) := f (a^0, \cdots, a^k | h)   \tr (A^0\cdots  A^k U(h)), \quad a^i\otimes A^i\in \maA\otimes \End(X)\text{ and } h\in H.
$$

\begin{lemma}
The cochain $f \sharp \tr$ is  continuous and equivariant, and we have
$$
b(f \sharp \tr) = (bf) \sharp \tr\text{ and } f \text{ cyclic }\Rightarrow f\sharp \tr\text{ cyclic}.
$$
\end{lemma}

\begin{proof}
Continuity is obvious. We have for $a^i\in \maA$, $A^i\in \End(X)$ and $h,g\in H$,
$$
(f \sharp \tr) (h (a^0 \otimes A^0), \cdots, h(a^n\otimes A^n) | h g h^{-1})
$$
\vspace{-0.7cm}
\begin{eqnarray*}
&=&(f \sharp \tr) (h a^0 \otimes U(h) A^0 U(h^{-1}), \cdots, ha^n\otimes U(h) A^n) U(h^{-1} | h g h^{-1})\\  & = & f(ha^0, \cdots, h a^n | h g h^{-1}) \tr [(U(h)A^0 U(h^{-1}) \cdots (U(h) A^n U(h^{-1}) U(hgh^{-1})]\\
& = & f(ha^0, \cdots, h a^n | h g h^{-1}) \tr [U(h) (A^0 \cdots A^n)  U(g h^{-1})]\\
& = & f(a^0, \cdots,  a^n |  g ) \tr [A^0 \cdots A^n  U(g)]\\
& = & (f\sharp \tr) (a^0 \otimes A^0, \cdots, a^n\otimes A^n | g ).
\end{eqnarray*}
If we assume in addition that $f$ is cyclic, then
\begin{eqnarray*}
(f\sharp \tr)(a^n\otimes A^n, a^0\otimes A^0, \cdots, a^{n-1}\otimes A^{n-1} | g ) & = & (-1)^n f(a^0, \cdots , g(a^n) | g) \tr(A^n A^0\cdots A^{n-1} U(g))\\
& = & (-1)^n f(a^0, \cdots , g(a^n) | g) \tr (A^0\cdots A^{n-1} U(g) A^n)\\
& = & (-1)^n f(a^0, \cdots , g(a^n) | g) \tr (A^0\cdots A^{n-1} g(A^n) U(g))\\
& = & (-1)^n (f\sharp \tr) (a^0\otimes A^0, \cdots, g(a^{n}\otimes A^{n}) | g ).
\end{eqnarray*}
The proof of the relation $b(f\sharp \tr)=bf\sharp \tr$ is similar and omitted.
\end{proof}

Recall that the $H$-equivariant $K$-theory of the $H$-algebra $\maA$ is defined in terms of $H$-invariant idempotents in $H$-algebras $\maA\otimes \End(X)$ for the action
$$
h (a\otimes A) := (h a)\otimes U(h)AU(h^{-1}).
$$
where $X$ runs through the finite dimensional unitary representations $U$ of $H$.  More presisely, we identify any two such $H$-invariant idempotents $e\in \maA\otimes \End(X)$ and $e'\in \maA\otimes \End(X')$, if there exists finite dimensional unitary representations $(V,Y)$ and $(V',Y')$ and $H$-invariant elements $x\in \maA\otimes \Hom(X\oplus Y, X'\oplus Y')$ and $y\in \maA\otimes \Hom(X'\oplus Y', X\oplus Y)$ such that 
$$
x y = e'  \oplus 0\in \maA\otimes \End (X'\oplus Y') \text{ and } y x = e \oplus 0 \in \maA\otimes \End (X\oplus Y).
$$
This is an equivalence relation and the quotient is a monoid for direct sums, whose associated Grothendieck group is by definition the $H$-equivariant $K_0$-theory or simply $K$-theory of the $H$-algebra $\maA$. We denote it as usual by $K_0^H (\maA)$ or simply $K^H(\maA)$. It is clear from the definition  that when $\maA=\C$ with the trivial action, $K^H(\maA)$ is isomorphic to the representation ring $R(H)$ of $H$, that is the Grothendieck group associated with the monoid of finite dimensional unitary representations of $H$.  When the algebra $\maA$ is not unital, we add a unit with a trivial action of $H$ and get a unital $H$-algebra ${\widetilde \maA}$, and the character $\ep: {\widetilde \maA}\to \C$ induces by a straightforward functoriality $\ep_*:K^H({\widetilde \maA}) \to R(H)$. The equivariant $K$-theory of the non unital algebra $\maA$ is then by definition the kernel of $\ep_*$.

The following is a generalization of Proposition 14 in \cite{C85} and we give a proof based on a straightforward extension of Connes method.  We point out that this proposition can also be deduced from the results of  \cite{Voigt}.

\begin{proposition}\label{Pairing} 
Let $f$ be a continuous equivariant cyclic $2k$-cocycle  on the pair $(\maA, H)$. 
For any  finite dimensional unitary representation $U: H \to \maU(X)$, and any $H$-invariant 
idempotent $e$ in ${\mathcal A}\otimes \End (X)$, we define $< f, e >$,  by the formula
\begin{Equation}\label{pairing1}
\hspace{3cm}       $< f, e > (h) := (f\sharp \tr)(e, \cdots , e | h)$
\end{Equation}
\noindent
Then 
\begin{itemize}
 \item $<f,e>$ is a continuous central function on the compact group $H$.
\item Formula \ref{pairing1} induces a pairing between the equivariant $K$-theory and the equivariant cyclic cohomology of the $H$-algebra ${\mathcal A}$, i.e.
$$
K_0^H(\maA) \otimes H^{ev}_\lambda (\maA, H) \longrightarrow C(H)^H,
$$
where $C(H)^H$ denotes the central continuous functions on $H$.
\end{itemize}
\end{proposition}

When $\maA$ is unital, the statement in the second item is clear. When $\maA$ is non unital, it needs some explanation.  As usual we denote by ${\widetilde{\maA}}$ the algebra $\maA \oplus \C$ endowed with the extended action (trivial on the $\C$ factor) and with the natural locally convex topology.  Then the second item means that if $X$ and $X'$ are two unitary representations of $H$ and if $\te = (e,\Lambda)$ and $\te'=(e',\Lambda')$ are two idempotents in ${\widetilde \maA} \otimes \End(X)$ and ${\widetilde \maA} \otimes \End(X')$, with $e\in \maA \otimes \End(X)$, $e'\in \maA \otimes \End(X')$, $\Lambda \in \End(X)$ and $\Lambda '\in \End(X')$, such that $[\te] - [\te']$ defines a class $x$ in $K^H(\maA)$, then the scalar
$$
(f\sharp \tr)(e, \cdots , e | h) - (f\sharp \tr)(e', \cdots , e' | h)
$$
only depends on the class $x$ and on the equivariant cyclic cohomology class of $f$. 

\begin{proof}
We give a brief proof for the convenience of the reader. We first notice that $f\sharp \tr$ is an equivariant continuous cyclic $2k$-cocycle on the $H$-algebra  ${\mathcal A}\otimes \End (X)$ and that $e$ is $H$-invariant. Therefore, we have
$$
(f\sharp \tr)(e, \cdots , e | g hg^{-1}) \,\,=\,\, (f\sharp \tr)(g^{-1}e, \cdots , g^{-1}e | h )  \,\,=\,\, (f\sharp \tr)(e, \cdots , e | h).
$$
Hence, $<f,e>(ghg^{-1})\,= \,\, <f,e>(h)$ and $<f,e>$ is central. The action of $H$ on $\maA$ being strongly continuous, the corresponding action on $ {\mathcal A}\otimes \End (X)$ is also strongly continuous. Now, continuity of $f$ implies continuity of $f\sharp \tr$. Therefore, for any $e$ as above, the map $<f,e>$ is a continuous function on $H$. 

In order to prove the second item, we extend  Connes' proof given in the non equivariant case in \cite{C85}. We can assume that $\maA$ is unital, see for instance \cite{BHII}. Using classical matrix techniques, we reduce the proof to showing that if $e,e'$ are $H$-equivariantly conjugated in $\maA$, then $<f,e>\,= \,\,<f,e'>$. We first prove that if $e'=ueu^{-1}$ with  $u$ an $H$-invariant element, then there exists an odd $H$-equivariant cyclic cochain $\varphi= \varphi_u$ such that  
\begin{Equation}\label{bphi}
	\hspace{2cm} $<f,e'> - <f,e>\,\,=\,\,<b\varphi,e>.$
\end{Equation}
For any $H$-invariant element $a\in {\mathcal A}$, the Hochschild cochain
$$
	f_a(a^0,...,a^{n-1} | h)=f (a^0,...,a^{n-1},a |h),
$$
is $H$-equivariant and we set
$$
	(\delta^*_a f)(a^0,...,a^n |h)=\sum_{0\leq i\leq n} f(a^0,...,aa^i - a^i a,...,a^n | h).
$$
Then $\delta^*_a f$ is an $H$-equivariant  cochain on $(\maA, H)$.  Using the relation $bf=0$, we easily get the following relation
$$
	(-1)^n (A \circ b')(f_a) (a^0, \cdots, a^n | h) = f((h^{-1} a) a^0 - a^0 a, a^1, \cdots, a^n |h) + \sum_{j=1}^n f(a^0, \cdots, a^{j-1}, aa^j - a^j a, a^{j+1}, \cdots, a^n)
$$
where $A$ is the cyclic operator on $C^*(\maA, H)$ defined in the previous section.
Since $a$ is supposed to be $H$-invariant, we deduce that $(-1)^n \delta^*_a f = (A \circ b')(f_a)$. Thus, using the relation $A\circ b'=b\circ A$, we obtain $\delta_a^* f=b [(-1)^n A (f_a)]$, which shows that the class of $\delta_a^* f$ in $H_\lambda^*(\maA, H)$ is trivial.

Now if $u$ is an invertible $H$-invariant element of ${\mathcal A}$, then  there exists an invertible $v\in M_2(\maA)$ and $b\in M_2(\C)$  such that $u=e^{vbv^{-1}}:= ve^b v^{-1}$, see \cite{C85} again.  We set $a=vbv^{-1}$ and write $\rho:M_2(\maA)\to M_2(\maA)$ for conjugation by $e^a$.   Then a straightforward computation shows that
$$
	\rho^*f-f \,\, = \,\,  (b\circ A)(\sum_{k\geq 1} 
	{\frac{(\delta_a^*)^{k-1}(f_a)}{k!}}),
$$
finishing the proof of relation \ref{bphi}. Notice that the convergence makes sense using only the convergence in $M_2(\C)$. Now, we have for any $h\in H$ and using that $e$ is $H$-invariant and $n$ is even
$$
<b\varphi,e> (h) \,\, = \,\, b'\varphi (e, \cdots, e|h) + (-1)^n \varphi( (h^{-1}e) e, e\cdots, e) \,\, = \,\,$$
$$
\varphi (e, \cdots, e | h) \sum_{j=0}^n (-1)^j
\,\, = \,\, \varphi (e, \cdots, e | h).
$$

Finally, using cyclicity of $\varphi$, we deduce that 
$$
\varphi (e, \cdots, e | h)\,\, = \,\, - \varphi (e, \cdots, e, h(e) | h) \,\, = \,\,- \varphi (e, \cdots, e, e | h).
$$
\end{proof}

\begin{remark}
The pairing in the odd case can  be defined similarly extending Proposition 15 in \cite{C85} to the equivariant setting.
\end{remark}

\begin{remark}\ The evaluation $<f,e>(h)$ of the pairing at $h$ is well defined provided the idempotent $e$ is only $h$-invariant. It does not depend on the group $H$ as far as this latter exists. 
\end{remark}

In the rest of this section, we explain the relation with the equivariant index pairing in the case of equivariant Fredholm modules.   In particular, an important example of an equivariant cyclic cocycle arises from index theory and is given as the equivariant Chern-Connes character of an equivariant finitely summable Fredholm module. We concentrate on the even case, as the odd case is similar.  Let $\maA$ be a locally convex $H$-algebra endowed with a continuous representation as operators on a $\Z_2$-graded Hilbert space $\maH=\maH^+\oplus \maH^-$. Let $\gamma$ be the grading involution, and assume that $\gamma$ commutes with the elements of $\maA$. Let the compact group $H$ act on $\maH$ by unitaries commuting with the grading $\gamma$ so that  $H$ acts on $\maA$ by conjugation. Let $F$ be an $H$-invariant symmetry ($F^*=F$ and $F^2=I$) on $\maH$ such that
$$
F\circ \gamma + \gamma \circ F = 0 \text{ and } [F, a] = F\circ a - a \circ F \in L^{2p} (\maH),
$$
where $L^{q} (\maH)$ is as usual the $q$-th Schatten ideal of operators $T$ such that $(T^*T)^{q/4}$ is a Hilbert-Schmidt operator. For $T\in L^1(\maH)$, we denote by $\tr(T)$ the trace of the operator $T$.

\begin{proposition}\label{ChernConnes}
The map
$$
(a^0 \otimes \cdots \otimes a^{2p} | h )\longmapsto (-1)^p \tr (\gamma a^0 [F,a^1]...[F,a^{2p}] U(h) ).
$$
is a continuous equivariant cyclic cocycle on the pair $(\maA, H)$.
\end{proposition}

\begin{proof}
Denote by $f$ this map and forget the $(-1)^p$. Continuity is obvious, given our assumptions. We compute 
\begin{eqnarray*}
f(ga^0, \cdots, ga^{2p} | ghg^{-1}) & = & \tr (\gamma  U(g) a^0 [F,a^1]\cdots [F,a^{2p}] U(g^{-1}) U(g) U(h) U(g^{-1}))\\
& = & \tr (U(g) \gamma  a^0 [F,a^1]\cdots [F,a^{2p}] U(h) U(g^{-1}))\\
& = & \tr (\gamma  a^0 [F,a^1]\cdots [F,a^{2p}] U(h) ).
\end{eqnarray*}
On the other hand, using the equality $\tr (\gamma [F, \omega]) =0$ for $\omega\in B(\maH)$ such that $\gamma \omega = (-1)^{\pa \omega} \omega \gamma$, the equivariant cyclicity is proved as follows:
\begin{eqnarray*}
f(a^{2p}, a^0, \cdots, a^{2p-1} | h) &=& \tr (\gamma a^{2p} [F,a^0]\cdots [F,a^{2p-1}] U(h) )\\
& = & \tr (\gamma  h(a^{2p}) U(h) [F,a^0]\cdots [F,a^{2p-1}])\\
& = & \tr (\gamma  [F, h(a^{2p}) U(h) a^0]\cdots [F,a^{2p-1}]) - \tr (\gamma  [F, h(a^{2p})] U(h) a^0 [F, a^1]\cdots [F,a^{2p-1}])\\ 
& = & - \tr (\gamma  [F, h(a^{2p})] U(h) a^0 [F, a^1]\cdots [F,a^{2p-1}])\\
& = & \tr ( [F, h(a^{2p})] U(h)\gamma  a^0 [F, a^1]\cdots [F,a^{2p-1}])\\
& = & \tr (\gamma  a^0 [F, a^1]\cdots  [F, h(a^{2p})] U(h))\\
& = & f(a^0,  \cdots, a^{2p-1}, h(a^{2p}) | h).
\end{eqnarray*}
To finish the proof of the proposition, it remains to show that $bf=0$. But, we have
\begin{eqnarray*}
(b' f) (a^0, \cdots, a^{2p+1} | h) & = & \tr (\gamma a^0 a^1 [F, a^2] \cdots [F, a^{2p+1}] U(h)) \\
& + & \sum_{j=1}^{2p} (-1)^j \tr (\gamma a^0 [F, a^1]\cdots ([F, a^j] a^{j+1} + a^j [F, a^{j+1}]) \cdots [F, a^{2p+1}] U(h))\\
& = & \tr (\gamma a^0 [F, a^1] \cdots [F, a^{2p}] a^{2p+1} U(h))\\
& = & \tr (\gamma a^0 [F, a^1] \cdots [F, a^{2p}] U(h) h^{-1} (a^{2p+1}) )\\
 & = & \tr (\gamma h^{-1} (a^{2p+1}) a^0 [F, a^1] \cdots [F, a^{2p}] U(h)  ),
\end{eqnarray*}
which completes the proof.
\end{proof}

\begin{definition}
We define the equivariant Chern-Connes character $\ch^H(\maH, F, \gamma)$ of $(\maH, F, \gamma)$ as the equivariant cyclic cohomology class of the equivariant cyclic cocycle defined in Proposition \ref{ChernConnes}. 
\end{definition}

Now let $({\mathcal H},F,\gamma)$ be 
an even $2p$-summable $H$-equivariant Fredholm module over the  $H$-algbera ${\mathcal A}$ as defined above.  For any finite dimensional unitary representation $U: H \to \maU(X)$ of the compact  group $H$, and any $H$-invariant projection
$e\in {\mathcal A}\otimes \End(X)$, the operator   $e\circ [F\otimes id_X]\circ e$,
acting on $e(H\times X)$, is an $H$-invariant Fredholm operator which anticommutes with the grading $\gamma$ of 
${\mathcal H}=\maH^+ \oplus \maH^-$ and induces the Fredholm operator 
$$
F_e^+ := e\circ [F\otimes id_X]\circ e : e(\maH^+ \otimes X) \longrightarrow e(\maH^- \otimes X).
$$
The $H$-equivariant index of $F_e^+$ is then an element of the representation ring $R(H)$ of  $H$. The map $e\mapsto \Ind (F_e^+)$ induces the map $\Ind^{H}_{{\mathcal H},F, \gamma}$,
$$
	\Ind^{H}_{{\mathcal H},F, \gamma}: K^{H}({\mathcal A}) \longrightarrow R(H).
$$

Recall that we have previously defined the equivariant Chern-Connes character $\ch^H(\maH, F, \gamma)$ of $(\maH, F, \gamma)$, as an equivariant  cyclic cohomology class over the $H$-algebra $\maA$.

\begin{proposition}\  We have
$$
	\Ind^{H}_{{\mathcal H},F, \gamma} (e) = <\ch^H(\maH, F, \gamma) , [e] >,
$$
where $[e]$ is the equivariant $K$-theory class of $e$.
\end{proposition}

\begin{proof}\ 
Let $X$ a finite dimensional representation of $H$ and let $e\in A\otimes \End(X)$ be an $H$-invariant idempotent. Then $e$ acts as an even degree operator on $\maH\otimes X$ and we  denote by $P$ the $H$-invariant Fredholm operator $F_e^+$ defined above. It is then easy to check that the operator $Q$ which is the same operator $e(F\otimes id_X)e$ but acting from $\maH^-\otimes X$ to $\maH^+\otimes X$,  is an $H$-invariant parametrix for $P$ modulo the Schatten ideal $L^{p}$.   Now recall the equivariant Atiyah-Bott formula for any $h\in H$, \cite{BenameurWorld},
$$
	 \Ind^{H}(P)(h) \,\, =\tr_{e({\mathcal H}^+\oplus X)}( (1-QP)^{p}\circ U(h))-\tr_{e({\mathcal H}^-\oplus X)}
	((1-PQ)^{p}\circ U(h)).
$$
Using the $H$-invariance of  $e$ and $F$ and the relation $
e [F\otimes id_X, e] e = 0,$
we deduce that 
$$
e-e(F\otimes id_X)e(F\otimes id_X)e = - e[F\otimes id_X, e ]^2.
$$
The computation of $[e-e(F\otimes id_X)e(F\otimes id_X)e]^p$ then finishes the proof. 
\end{proof}

As a consequence of the previous proposition, we deduce the following integrality result.

\begin{corollary}\ Let $({\mathcal H},F, \gamma)$ be a finitely summable $H$-equivariant Fredholm module over $\maA$ as above.  Then 
$$
	\left<\ch^H(\maH, F, \gamma), K^{H}({\mathcal A})\right>\,\, \subset \,\, R(H).
$$
\end{corollary}

When $H$ is the finite cyclic group of order $n$, we get
$$
	\left<\ch^H(\maH, F, \gamma), K^{H}({\mathcal A})\right> (h)  \,\, \subset \,\, \Z[e^{2i\pi/n}], \quad \forall h\in H.
$$
Therefore, for foliated involutions, one gets the integrality result $\left<\ch^H(\maH, F, \gamma), K^{H}({\mathcal A})\right> \subset \Z$.

\section{Haefliger currents and leafwise diffeomorphisms}

We now focus on foliations. We recall the correspondence between Haefliger homology and Connes' periodic cyclic cohomology, and then prove that the cyclic cocycles arising from this correspondence are equivariant with respect to leafwise actions of leafwise volume preserving diffeomorphisms.

As above, $(V,F)$ is a smooth compact foliated manifold of codimension $q$, 
and ${\mathcal G}$ is the holonomy groupoid of $F$.   Denote the metrics on $\cG_x$ and $\cG^y$ by $dvol_r = r^*(dvol_F)$  and $dvol_s = s^*(dvol_F)$, where $dvol_F$ is the given metric on the leaves of $F$.  Recall that  $\nu \subset TV$ is the normal bundle to $TF$,  
so $TV = TF \oplus \nu $. Any leafwise path $\gamma:[0,1] \to V$ defines a germ of diffeomorphisms denoted $h_{\gamma}$ (and called the holonomy germ) from transversals at $\gamma(0) = s(\gamma)$ to transversals at $\gamma(1) = r(\gamma)$.  The differential of the germ, also denoted $h_{\gamma}$, gives a well defined linear isomorphism $h_{\gamma}:\nu_{s(\gamma)} \to \nu_{r(\gamma)}$. This defines an action of $\cG$ on $\nu$.  Transposing this action gives an action of $\cG$ on the dual normal bundle $\nu^*$. Notice that $\nu^*$ is naturally defined as the subbundle of $T^*V$
$$
\nu^* \,\,=\,\, \{\alpha\in T^*V \,\,\, | \,\,\, \forall X\in TF, \,\,\, \alpha(X) \, = \, 0\}.
$$
Thus this action of $\cG$ on $\nu^*$, as well as its action on all the exterior powers $\wedge^k \nu^* $ for $0\leq k\leq q$, is natural. 

Denote by $\Omega^j$ the space $C_c^{\infty}(\cG,\wedge^j r^*(\nu^*))$.  This is a graded algebra where
\begin{Equation}\label{prec}
$\dd   \hspace{1cm}
\forall (\omega_1,\omega_2)\in \Omega^j\times \Omega^{j'} \,\, \text{and } \,\,\forall \gamma\in \cG^y, 
\quad \quad 
\omega_1\omega_2(\gamma)  \,\,= \,\,
\int_{\cG^y} \omega_1(\gamma_1)\wedge \wh_{\gamma_1}(\omega_2(\gamma_1^{-1}\gamma))dvol_s.
$
\end{Equation}

\noindent
The push-forward map $\wh_{\gamma_1}:  \wedge^{j'} \nu^*_{s(\gamma_1)}  \to \wedge^{j'}  \nu^*_{r(\gamma_1)}$  is $\wh_{\gamma_1} = [h_{\gamma_1}^*]^{-1}   =  h_{\gamma_1^{-1} }^*$.
{
There is a transverse differential \cite{C94} p.\ 266,
$$
	d_\nu : \Omega^j \to \Omega^{j+1},
$$
which satisfies all the properties of a precycle (definition below) required of it, which we now recall.

For $\gamma\in \cG$ and $X_1, \cdots , X_j\in \nu_{r(\gamma)}$, there exist  unique tangent vectors $Y_1, \cdots , Y_j\in \nu_{\cG, \gamma}$ such that
$$
r_* Y_i = X_i \text{ and } s_* Y_i = h_{\gamma^{-1}} X_i.
$$
If $d_\cG$ is the de Rham differential on the smooth manifold $\cG$, then for $\omega \in \Omega^j$, 
$$
d_\nu  \omega (X_1, \cdots, X_j) := d_\cG \omega (Y_1, \cdots, Y_j).
$$
In  \cite{C94} p.\ 267, Connes showed that  $d_\nu $ is a graded differential of degree $1$ on $\Omega=\oplus \Omega^j$, which satisfies the relation
$$
d_\nu ^2 (\omega) = [\theta, \omega],
$$
with $\theta$ a generalized (compactly supported) section of $r^*(\wedge^2\nu^*)$ over $\cG$, which is a multiplier of $\Omega$, i.e. for any $\omega \in \Omega$, $\theta\omega$ and $\omega\theta$ make sense in $\Omega$. In addition, $d_\nu  (\theta)$ is zero.

Fix a finite distinguished open cover $\{U_i \}$ for the foliation and choose smooth transversals $T_i\subset U_i$ such that ${\overline {T_i}} \cap {\overline {T_{i'}}} = \emptyset$ for $i\not=i'$. Then with $T= \bigcup_i T_i$, the space $\Omega_c(T/F)$ of Haefliger forms on the complete transversal $T$ is the quotient of the space of compactly supported differential forms on the smooth manifold $T$ by the closure of the subspace generated by all forms, when well defined,  of the form $h_\gamma^*\alpha - \alpha$.  See \cite{Ha}. 
A holonomy invariant $k$-current $C$ assigns a real number to any compactly supported differential $k$ form defined on any transversal, with the stipulation that $C(h_\gamma^*\alpha - \alpha)=0$.
Any such $C$ gives a continuous (for the smooth topology) linear form on $\Omega_c^k(T/F)$, and such a form is called a Haefliger current.  Let $dvol_F$ be the leafwise volume form determined by the restriction of a metric $g$ on $V$ to the leaves of $F$. Any $\omega\in \Omega^k$ can be restricted to the units $\cG^{(0)}=V$ to yield a smooth differential form $\omega|_V$ on $V$, which is actually a section of $\wedge^k \nu^*$. Given $\alpha \in C^{\infty}(V;\wedge^k \nu^*)$, integration over the leaves of $F$ of  $\alpha \wedge dvol_F$, as defined for instance in \cite{Ha}, yields
$$
\int_F: C^\infty(V; \wedge^k\nu^*) \longrightarrow \Omega^k_c(T/F).
$$
Set
$$
\int_C \;\; \omega := \left<\int_F (\omega |_V) , C \right>.
$$

Recall from \cite{C85} the definition of a precycle over an algebra $A$: 
\begin{definition}\label{precyc} Let $A$ be a $\C$-algebra. A $k$-precycle over $A$ is a quadruple $(\Omega,d,\int,\theta)$ such that:
\begin{enumerate}
\item $\Omega$ is a graded algebra $\Omega=\Omega^0\oplus ...\oplus \Omega^k$,  the product
sends $\Omega^j\otimes \Omega^{j'}$ into $\Omega^{j+j'}$, and $\Omega^0=A$;
\item $d$ is a graded differential of degree 1 on $\Omega$, so $d:\Omega^j \to \Omega^{j+1}$ and
$$
	d(\omega_1\omega_2)=(d\omega_1)\omega_2+(-1)^{j}\omega_1(d\omega_2) \quad\forall (\omega_1,\omega_2)\in \Omega^j\times \Omega,
$$
\item  $\theta\in \Omega^2$ (see remark \ref{theta} below) satisfies $d\theta = 0$, and for all $\omega$,
$$
	d^2(\omega)=[\theta,\omega]= \theta \omega - \omega \theta;
$$
\item $\dd \int:\Omega^k\to \C$ is a closed graded trace, i.\ e.\
$$
	\forall (\omega_1,\omega_2)\in \Omega^j\times \Omega^{k-j}, \quad
	\int \omega_1\omega_2=(-1)^{j(k-j)}\int \omega_2\omega_1; \quad 
	\text{ and } \,\, \forall \omega\in \Omega^{k-1},  \quad \int d\omega =0.
$$
\end{enumerate}
A $k$-precycle over $A$ is a $k$-cycle if  $d^2=0$.
\end{definition}
\begin{remark}
 The condition $\Omega^0=A$ is not necessary, even though it will be sufficient for our applications, and one can replace it by the existence of an algebra morphism from $A$ to $\Omega^0$.
\end{remark}
\begin{remark} \label{theta}The condition $\theta\in \Omega^2$ is too strong and Connes' $X$-trick (see below) works just as well  even if $\theta$ is only a  degree $2$ multiplier of $\Omega$. In general, it suffices to assume the existence of  $\theta$ such that $d\theta$ and $[\theta,.]$ do make sense. 
With this in mind, we proved in \cite{BHI} that $(\Omega, d_\nu , \int_C, \theta)$ is a $k$-precycle over the convolution algebra $C_c^\infty (\cG)$ of smooth compactly supported functions on the graph $\cG$. 
\end{remark}
The classical example is $A=C_c^{\infty}(M)$, where $M$ is  a smooth $m$ dimensional manifold.   Then every closed de Rham $k$-current $C$ on $M$ gives rise to a $k$-cycle over $A$, by considering the compactly supported de Rham complex, truncated at the level $k$, together with the closed graded trace induced by $C$.
Other examples come from the study of finitely summable Fredholm modules over $\C$-algebras, see  \cite{C85} for the details.

Connes has given a curvature method, nowadays called Connes' $X$-trick, \cite{C94} p.\ 229, which assigns to any  $k$-precycle $(\Omega, d, \int, \theta)$ over an algebra $A$ a $k$-cycle $({\widetilde \Omega} , \delta, \Phi, 0)$ over an extension of $A$ and hence over $A$. Since the $X$-trick will be important in the sequel, we recall it for the convenience of the reader.  As a vector space, the new graded algebra ${\widetilde \Omega} = M_2(\Omega)$, the $2$ by $2$ matrices over $\Omega$, but the product is not the usual one.  Set
$$
      \Theta=\left(
      \begin{array}{cc}
      1 &\quad 0\\
      0 & \quad \theta \end{array} \right)
$$
and define the product $*$ on ${\widetilde \Omega}$ by setting 
$$
      {T} * {T}' := {T} \Theta {T}'.
$$
Denote by $\pa \omega$ the degree of a homogeneous element $\omega \in \Omega$, and recall that $0$ has all degrees.  An element 
${T} \in {\widetilde \Omega}$ is homogeneous of degree $\pa 
{T} = k$ if
$$
      k =  \pa {T}_{11} = \pa {T}_{12}+1 =
      \pa {T}_{21} +1 = \pa  {T}_{22} +2.
$$
The differential $\delta$ on ${\widetilde \Omega}$ is defined 
on homogeneous elements of ${\widetilde \Omega} $ as
$$
\delta {T} = \left(\begin{array}{cc}
d {T}_{11} & d{T}_{12}\\
 - d{T}_{21} & - d {T}_{22}
\end{array}\right) 
+ 
\left(\begin{array}{cc} 0 & -\theta\\
1 & 0\end{array} \right){T} 
+ 
(-1)^{\pa{T}} {T}
\left( \begin{array}{cc}
 0 & 1\\ -\theta & 0
\end{array}\right).
$$
A straightforward computation then shows that $\delta ^2=0$. This makes $({\widetilde \Omega} , \delta)$ into a graded differential algebra.

The graded algebra $\Omega$ embeds as a subalgebra of 
${\widetilde \Omega} $ by using the homogeneous map
$$
	\omega \hookrightarrow \left( \begin{array}{cc} \omega & 0 \\ 0 & 0
	\end{array} \right).
$$ 

For homogeneous ${T}\in {\widetilde \Omega}^k$ define
$$
	\Phi (T) = \int {T}_{11} - (-1)^{\pa {
	T}} \int {T}_{22} \theta,
$$
and extend to arbitrary elements of ${\widetilde \Omega}^k$ by linearity. Then the linear form $\Phi$ is a closed graded trace on $({\widetilde \Omega} , \delta)$, see  \cite{C94}. Thus $({\widetilde \Omega} , \delta, \Phi, 0)$, denoted $({\widetilde \Omega} , \delta, \Phi)$ for short,  is a $k$-cycle. The precycle  $(\Omega, d_\nu, \int_C, \theta)$ defined in \cite{BHI} and recalled in Remark \ref{theta}, gives rise to the cycle that we denote $({\widetilde\Omega}, \delta_\nu, \Phi_C)$.

Given a general  $k$-cycle $(\Omega, d, \int)$ over an algebra $A$, define a cyclic $k$-cocycle (the Connes character of the cycle) by setting
\begin{Equation}\
\hspace{2in} 	$\dd \tau(a^0,...,a^k)=\int a^0 da^1 \cdots da^k.$
\end{Equation}

\begin{theorem}\label{tauC} \cite{BHI} Let $C$ be a closed Haefliger $k$-current for $(V,F)$ as above. Then the formula
$$
\tau_C(k_0,...,k_k)  \,\,:= \,\, 
\int_C \Bigl{(}[k_0 * \delta_\nu  k_1*\cdots *  \delta_\nu  k_k]_{11}|_V\Bigr{)} 
$$
defines a cyclic $k$-cocycle over the algebra $C_c^{\infty}(\cG)$, which is continuous with respect to the $C^\infty$ compact open topology.
\end{theorem}

\begin{proof}
 In \cite{BHI}, we proved that for any closed Haefliger $k$-current $C$, the quadruple $(\Omega, d_\nu , \int_C, \theta)$ defined above, is a precycle over the  convolution algebra of compactly supported smooth functions on $\cG$.
It is easy to check that for the associated cycle $({\widetilde \Omega},  \delta_\nu , \Phi_C)$, the multilinear functional $\tau_C$ is precisely the cyclic cocycle which is the Chern-Connes character of $({\widetilde \Omega},  \delta_\nu , \Phi_C)$.  Note that $\tau_C$ does not depend on the choice of $\nu$.

It remains to show continuity. But, the map $k\mapsto d_\nu k$ is clearly continuous with respect to the $C^\infty$ compact open topologies on $C_c^\infty (\cG)$ and $C_c^\infty (\cG, r^*\nu^*)$.  Multiplication by $\theta$ being continuous, $k\mapsto \delta_\nu k$ is also continuous for the induced topology on $M_2 ( C_c^\infty (\cG, r^*\nu^*))$. As we are only dealing with compactly supported sections, the $*$-product is continuous. Therefore, the multlinear map 
$$
(k_0, \cdots, k_r) \longmapsto [k_0 *\delta_\nu k_1 \cdots *\delta_\nu k_r]_{11}
$$
is continuous with respect to the $C^\infty$ compact open topologies on $C_c^\infty (\cG)$ and $C_c^\infty (\cG, r^*\wedge^*\nu^*)$. If we fix a smooth complete transversal $T$ to the foliation, then Haefliger integration over the leaves is known to be continuous from $\Omega^*(V)$ to $\Omega_c^* (T)$, endowed with the compact open topologies. By definition of a Haefliger current, we know that its action on $\Omega_c(T)$ is continuous. To sum up, we get continuity of the multilinear form
$$
(k_0, \cdots, k_r) \longmapsto \left< C, \int_F [k_0 *\delta_\nu k_1 \cdots *\delta_\nu k_r]_{11}|_V \wedge dvol_F \right>,
$$
which finishes the proof.
\end{proof}

\begin{remark}
For $k \geq 2$, $[k_0 * \delta_\nu  k_1*\cdots *  \delta_\nu  k_k]_{11}$ does involve the curvature $\theta$ of  the precycle $(\Omega, d_\nu , \int_C, \theta)$. 
\end{remark}
\begin{remark}
In the presence of an Hermitian bundle $E \to V$, it is easy to extend this construction, and to associate with any closed Haefliger current $C$ on the foliation $(V, F)$, a cyclic cocycle on the algebra $C_c^\infty (\cG, E)$ of smooth compactly supported sections over $\cG$ of the vector bundle  whose fiber at $\gamma\in \cG$ is as usual $\Hom(E_{r(\gamma)}, E_{s(\gamma)})$. See \cite{BHI} for the precise definition. 
\end{remark}

\medskip

Fix a diffeomorphism $f: V\to V$ which preserves the leaves and assume that $f$ preserves the leafwise Lebesgue measure defined by the volume form $dvol_F$.  This assumption will be referred to as the leafwise SL assumption. We proceed now to prove that the cyclic cocycle associated with a closed Haefliger current as above is indeed $f$-equivariant.  Our assumption on $f$  is satisfied for instance when $f$ is an element of a compact Lie group which acts on $V$ by leaf preserving diffeomorphisms, but this is not needed to prove the $f$-equivariance.  It is of course easy to construct examples on say $T^2$ where $f$ is not an element of a compact group and satisfies our leafwise SL assumption. 

We will need the following.
\begin{definition}\label{smosec}
The diffeomorphism $f$ is a holonomy diffeomorphism if there exists a smooth map $\varphi^f: V \to \cG$, so that for any $x\in V$, $s(\varphi^f (x)) = x$, $r(\varphi^f (x)) = f(x)$, and the holomony along $\varphi^f(x)$ coincides with the action of $f$ on transversals.
\end{definition}

\begin{lemma}
The diffeomorphism $f$ is a holonomy diffeomorphism in the following cases:
\begin{enumerate}
\item
When the holonomy is trivial, and the foliation is tame.
\item
When the foliation is Riemannian.
\item
When $f$ belongs to a connected Lie group which acts on $V$ by leaf-preserving diffeomorphisms.  More generally, if $f$ belongs to the path connected component of a holonomy diffeomorphism $g$ (for instance $g =$ identity) in the group of leaf-preserving diffeomorphisms.
\item
When restricted to the saturation $sat(V^f)$ of the fixed point submanifold $V^f$, that is the union of the leaves that intersect $V^f$. 
\end{enumerate}

\end{lemma}

\begin{proof}
For the first two items, see \cite{Hector2010}.    For the third, choose any smooth path $h_t$ from the $g$  to $f$ in the group, and define $\varphi^f(x)$ to be the element of $\cG$ determined the composition of $\varphi^g(x)$ with the leafwise path $t \to  h_t(x)$. 
For the last item, if $x\in sat(V^f)$, $y \in V^{f} \cap L_{x}$, and if $\gamma \in \cG^{y}_{x}$, set $\varphi^{f}(x) = (f \gamma^{-1})\circ \gamma$.  Then $\varphi^f$ is well defined, smooth,  depends only on $x$, and it works. 
\end{proof}

\begin{remark}
Theorem \ref{Lefschetz} shows that the Lefschetz class lives over $sat(V^f)$.  
\end{remark}

\begin{remark}
 A particular case of the second item  in the above lemma is when $f$ is given by the time one flow of a vector field on $V$ which is tangent to the foliation $F$. 
\end{remark}

To sum up, we are assuming from now on that $f$ is  a holonomy diffeomorphism which preserves the leafwise Lebesgue measure associated with $dvol_F$. Note that depending on which current we are using to produce Lefschetz formulae, the holonomy assumption will sometimes be unnecessary.  In particular, if the current is diffuse (in the sense that it is determined by its restriction to open dense subsets of a complete transversal) then we recall the
following classical density theorem due to Gilbert Hector, \cite{Hector}.

\begin{theorem}\label{density}
Let $T_1$ and $T_2$ be transversals of a foliation $F$, and consider the map between transversals $f:T_1 \to f(T_1) \subset T_2$.
Then there is a countable family $\{U_j\}$ of open subsets of $T_1$ so that:
\begin{enumerate}
\item $\dd \bigcup_j U_j$ is dense in $T_1$,
\item any restriction of $f$ to $f_j:U_j \to f(U_j)$ is a holonomy map, that is an element of the holonomy pseudo-group of $F$.
\end{enumerate}
\end{theorem}
\begin{proof}\  [Hector]
We may assume that $T_1$ and $T_2$ are subsets
of a global transversal $T$ and we denote by $\cP$ the pseudo-group of holonomy transformations of
$T$ induced by the foliation $F$. For any element $h \in \cP$, there exists a subset $U_h \subset T$, maybe the empty set, such that
$f(x) = h(x)$ for any $x\in U_h$ and this set is a closed subset of $T$ by continuity.
On the other hand for any $x\in T$ there exits an element $h_x \in \cP$ such that
$f(x) = h_x(x)$,  so
$\bigcup U_h$ covers T. But the set of maximal elements of $\cP$ is
countable, and therefore by Baire category theory there exists a set $U_{h_1}$ with non empty
interior $U_1$.   The result follows by a standard argument.
\end{proof}

Notice that in view of the Lefschetz problem we are interested in, $f$ will be an isometry of the ambiant Riemannian manifold $V$, in which case the SL assumption is of course satisfied.  The holonomy assumption is satisfied in all known examples. 
We now fix some Hermitian vector bundle $E$ over $V$ and assume that we have a bundle isomorphism
$$
A_f: f^*E \rightarrow E.
$$
We denote by  $\epsilon_{V,E}$  the Hilbert $C^*$-module over Connes' $C^*$-algebra $C^*(V,F)$ of the foliation, as defined in \cite{CS}. This is the Hilbert module associated with the $\cG$-equivariant continous field of Hilbert spaces $L^{2}(\cG_{x}, r^* E)$ over $V$, see \cite{C82}.
Define the endomorphism  
$\Psi^{E}(f) = \left(\Psi^{E}(f)_{x}\right)_{x \in V}$ of the Hilbert module $\epsilon_{V,E}$,
$\Psi^{E}(f)_{x} : L^{2}(\cG_{x},r^* E) \rightarrow L^{2}(\cG_{x}, r^* E)$,  
by setting
$$
\Psi^{E}(f)_{x}(\xi)(\gamma) \,\,=\,\, A^{-1}_f \Bigl{[} \xi( f^{-1}\gamma \cdot \varphi^{f^{-1}}(s(\gamma)))\Bigr{]} \,\, = \,\,
A^{-1}_f \Bigl{[} \xi(  \varphi^{f^{-1}}(r(\gamma)) \cdot \gamma)\Bigr{]},
$$
for $\xi \in L^{2}(\cG_{x}, r^*E)$ and $\gamma \in \cG_{x}.$  The second equality follows from the fact that the two paths $f^{-1} \gamma  \cdot   \varphi^{f^{-1}} (s(\gamma))$ and $\varphi^{f^{-1}} (r(\gamma))\cdot \gamma$ start and end at the same points and, since $f$ is a holonomy diffeomorphism,  the holonomy along them is the same.

In the special  case $E=V\times \C$ where $A_f$ acts as the identity of $\C$,   $\Psi^{E}(f)=\Psi(f)$ is, as we shall see soon, a multiplier of the $C^*$-algebra of the foliation which preserves $C_c^{\infty}(\cG)$.   
Recall that $k \in C_c^{\infty}(\cG)$ acts on $\oplus_{x\in V} L^2(\cG_x)$ as follows.  If $\xi  \in  L^2(\cG_x)$, and $\gamma \in \cG_x$, then
$$
k(\xi)(\gamma) \,\, = \,\, \int_{\gamma_1\in \cG_x}  k(\gamma \gamma_1^{-1})\xi( \gamma_1) dvol_r(\gamma_1).
$$
With this in mind, the operator $\Psi(f)\circ k$ has a smooth compactly supported kernel denoted by $k^{f,l}$. In the same way, $k \circ\Psi(f)$ has a smooth compactly supported kernel $k^{f,r}$. A simple computation gives 
$$
k^{f,l}(\gamma)    \,\, = \,\,       k(f^{-1} \gamma  \cdot   \varphi^{f^{-1}} (s(\gamma)))
 \,\, = \,\,       k(\varphi^{f^{-1}} (r(\gamma))\cdot \gamma), 
$$
where the second equality follows just as it does above.
Similarly,
$$k^{f,r}(\gamma) \,\, = \,\, k(\gamma \cdot \varphi^{f^{-1}} (f s(\gamma)))
 \,\, = \,\, k(  \varphi^{f^{-1}} (f r(\gamma))  \cdot f \gamma ).
$$
To see this note that 
$$
[k\circ \Psi(f)](\xi) \gamma \, \, = \,\, \int_{\cG_x} k(\gamma \gamma_1^{-1}) \xi( f^{-1}\gamma_1 \varphi^{f^{-1}} (x)) d\nu^x(\gamma_1).
$$
Set $\gamma'_1:= f^{-1}\gamma_1 \cdot \varphi^{f^{-1}} (x) \in \cG_x$. Then 
$$
\gamma_1^{-1} =  {\varphi^f(x)}^{-1} \cdot f{\gamma'}_1^{-1}= \varphi^{f^{-1}} (fx) \cdot f{\gamma'}_1^{-1} = {\gamma'}_1^{-1} \cdot \varphi^{f^{-1}} (f r(\gamma'_1))  = {\gamma'}_1^{-1} \cdot \varphi^{f^{-1}} (f s(\gamma {\gamma'_1}^{-1})).
$$
So
$$
[k\circ \Psi(f)](\xi) \gamma \, \, = \,\, \int_{\cG_x} k(\gamma {\gamma'}_1^{-1} \cdot \varphi^{f^{-1}} (f s(\gamma{\gamma'_1}^{-1}))) \xi(\gamma'_1)d \nu^x(\gamma'_1).
$$
As $[k\circ \Psi(f)](\xi) \gamma$ is also given by
$$
[k\circ \Psi(f)](\xi) \gamma \, \, = \,\, \int_{\cG_x} k^{f,r}(\gamma {\gamma'}_1^{-1}) \xi(\gamma'_1)d \nu^x(\gamma'_1),
$$
we have the formula.

Recall that the action of $f$ on $k \in C_c^{\infty}(\cG)$ is given by 
$$
(f k)(\gamma) \,\,=\,\,k(f^{-1}\gamma), 
$$  
where $\gamma \in \cG$.  Then we immediately have
\begin{lemma}\label{Haction} 
\hspace{2cm} $\Psi(f) \circ k \circ \Psi(f^{-1}) = f(k).$
\end{lemma}
\begin{proof}
For $\gamma\in \cG_x$,
$$
(k^{f^{-1},r})^{f,l} (\gamma) = k^{f^{-1},r} (f^{-1} \gamma  \circ   \varphi^{f^{-1}} (x)) = k (f^{-1} \gamma  \circ   \varphi^{f^{-1}} (x) \circ \varphi^{f} (f^{-1} x)) = k (f^{-1} \gamma).
$$
\end{proof}

It is then easy to check that the same relations hold in the presence of the general Hermitian bundle $E$ over $V$ with isomorphism $A_f$ as above. So for instance, we have for any $k\in C_c^\infty (\cG, E):=C_c^\infty (\cG, \Hom(s^*E, r^*E))$:
$$
\Psi^E(f) \circ k \circ \Psi^E(f^{-1}) = f(k), \text{ where } f(k) (\gamma)= (A_f)_{r(\gamma)}^{-1} \circ k(f^{-1}\gamma) \circ (A_f)_{f^{-1}s(\gamma)}.
$$

\begin{proposition}\label{ccoc} Let $C$ be a closed Haefliger current and  $\tau_C$ the cyclic cocycle on $C_c^{\infty}(\cG, E)$ associated with $C$ by Theorem \ref{tauC}.   Then $\tau_C$ satisfies the following equation.
$$
\tau_C(k_0, ...,k_i,\Psi^E(f) \circ k_{i+1},  ..., k_r)   \,\, = \,\, 
\tau_C(k_0, ...,k_i \circ \Psi^E(f), k_{i+1}, ... , k_r).
$$
\end{proposition}

\begin{proof}\ 
We shall give the proof for the trivial line bundle $E=V\times \C$, from which the general case is  easily deduced.   If $C$ is a $0$ current, we must show $\tau_C(\Psi(f) \circ k) =
\tau_C(k\circ \Psi(f))$.  But $k\circ \Psi(f) =  \Psi(f^{-1}) \circ  \Psi(f)\circ k\circ \Psi(f) = f^{-1}(\Psi(f)\circ k)$, and in the proof of Proposition \ref{ccoc2}, we show that $\tau_C(f^{-1}(\Psi(f)\circ k))  = \tau_C(\Psi(f)\circ k)$. 

Now suppose that $r=\dim (C)>0$.  Since $\tau_C$ is cyclic, we need only show
$$
\tau_C(k_0,  \Psi(f)  \circ k_1, k_2, ..., k_r)   \,\, = \,\, 
\tau_C(k_0 \circ \Psi(f), k_1, ... , k_r).
$$ 
Recall that the product $k_0 *  \delta_\nu  k_1 *  \delta_\nu  k_2 * \cdots *  \delta_\nu  k_r$  in Connes' $X$-trick 
is given by
$$
k_0 *  \delta_\nu  k_1 *  \delta_\nu  k_2 * \cdots *  \delta_\nu  k_r \,\, = \,\, 
\left( \begin{array}{cc} k_0 & 0 \\ 0 & 0  \end{array} \right)  
\left(\begin{array}{cc}d_\nu  k_1 & k_1\\ \theta k_1 & 0 \end{array}\right)
\left(\begin{array}{cc}d_\nu  k_2 & k_2\\  \theta k_2 & 0 \end{array}\right)
\cdots
\left(\begin{array}{cc}d_\nu  k_r & k_r\\  \theta k_r & 0 \end{array}\right),
$$
where $\theta = d_\nu^2$, the curvature of $d_{\nu}$ defined above.
To compute $\tau_C(k_0,k_1, ..., k_r)$, we need only the $(1,1)$ component of this product.

There are two types of terms.  The first type consists of terms of the form 
$k_0 \, (\Psi(f)\circ k_1) \, A$, where $A$ is a polynomial in the variables $k_2,  \ldots, k_r,  d_\nu  k_2,  \ldots,  d_\nu  k_r,$ and $\theta$.    It is immediate that $k_0 \, (\Psi(f)\circ k_1) = (k_0  \circ \Psi(f))\,k_1 $, that is $k_0 k_1^{f, \ell} = k_0^{f, r} k_1$, so we have the result for these terms. 
The second type consists of terms of the form 
$k_0 \, d_\nu  (\Psi(f)\circ k_1) \, A =  k_0 \, d_\nu  k_1^{f, \ell} \, A $.

\begin{lemma}
For $k \in C_c^{\infty}(\cG)$,
$$
d_\nu  k^{f,\ell} (\gamma) \,\, = \,\, ^t(f^{-1}_*) ( d_\nu  k (\varphi^{f^{-1}}(r(\gamma) \circ \gamma)))
\text{\quad and \quad} 
d_\nu  k^{f,r} (\gamma)   \,\, = \,\,
d_\nu k (\gamma\circ \varphi^{f^{-1}}(fs(\gamma))).
$$ 
\end{lemma}

\begin{proof}

Let $\gamma \in \cG^y_x$,  and set 
$\gamma' =\varphi^{f^{-1}}(r(\gamma)) \circ \gamma \in  \cG^{f^{-1}(y)}_x$.
Let $k \in C_c^{\infty}(\cG)$.
Note that 
$d_\nu  k^{f,l}(\gamma) \in r^*(\nu^*)_{\gamma} = \nu^*_y$ and 
$d_\nu  k(\gamma') \in r^*(\nu^*)_{\gamma'} = \nu^*_{f^{-1}(y)}$.
Given $X \in \nu_y = r^*(\nu)_{\gamma}$, let $Y \in T\cG_{\gamma}$ be the unique  vector with $r_*(Y) = X$ and $ s_*(Y) = h_{\gamma^{-1}}X$.  Then 
$$
d_\nu  k^{f,l}(X)\,\,=\,\,  d_{\cG}k^{f,l}(Y).
$$

The vector $X$ also gives the vector $Y' \in T\cG_{\gamma'}$, determined by the requirements that  
$r_*(Y') =  X'$ where $X' = f^{-1}_*(X) = h_{\varphi^{f^{-1}} (r(\gamma))}X$, and $ s_*(Y') = h_{\gamma^{-1}}X = h_{(\varphi^{f^{-1}}(r(\gamma)) \circ \gamma)^{-1}}X' = h_{ \gamma^{'-1}}X'$. 
Then
$$
d_\nu  k(X')\,\,=\,\,  d_{\cG}k(Y').
$$
Let $\gamma_t:[0,1]\to V$ for $t$ near $0$ be a smooth family of leafwise paths which defines the tangent vector $Y$.  It is sufficient for  $\gamma_t$ to satisfy three requirements:  
$$\gamma_0 = \gamma; \quad \frac{d\gamma_t}{dt}(0) \, |_{t = 0}  = h_{\gamma^{-1}}X;  \quad   \frac{d\gamma_t}{dt}(1) \, |_{t = 0} = X.
$$
Now consider the smooth family of leafwise paths  $\gamma'_t =  \varphi^{f^{-1}} (r(\gamma_t))\circ \gamma_t$.  First $\gamma'_0 = \gamma'$.   Second, $s(\gamma'_t) = s(\gamma_t)$ so $\dd  \frac{d\gamma'_t}{dt}(0) \, |_{t = 0} =  \frac{d\gamma_t}{dt}(0) \, |_{t = 0} = h_{\gamma^{-1}}X  = h_{\gamma^{'-1}}X'$. Finally, $r(\gamma'_t) =f^{-1}(r(\gamma_t))$, so $\dd \frac{d\gamma'_t}{dt}(1) \, |_{t = 0}  = f^{-1}_*(X) = X'$.
So the family  $\gamma'_t$ defines the vector $Y'$.  
As $k^{f,l}(\gamma_t) \, =  \,  k(\gamma'_t)$, we have
$$
d_\nu  k^{f,l}(X)\,=\, d_{\cG}k^{f,l}(Y) \,=\, 
\frac{d(k^{f,l}(\gamma_t))}{dt}  \, | _{t=0} \, =\,  
\frac{d(k(\gamma'_t))}{dt} \, |_{t=0} \,=\,
$$
$$ d_{\cG}k(Y')  \,=\,  d_\nu  k(X') \,=\,  d_\nu  k( f^{-1}_*(X))
\,=\,  (^t(f^{-1}_*) d_\nu  k)(X).
$$
Therefore, for all $\gamma \in \cG$,
$$
d_\nu  k^{f,l}(\gamma)   \,\,=\,\,  ^t(f^{-1}_*)(d_\nu (k)(\gamma'))    \,\,=\,\,
^t(f^{-1}_*) ( d_\nu  k (\varphi^{f^{-1}}(r(\gamma) \circ \gamma))).
$$

A similar argument proves the second relation. 
\end{proof}

To finish the proof of Proposition \ref{ccoc}, we have.
\begin{lemma}\label{switchd}
For $k_0,k_1 \in C_c^{\infty}(\cG)$,
$$
k_0 d_\nu k^{f,l}_1 \,\, = \,\,  k_0^{f,r} d_\nu k_1.
$$
\end{lemma}
	
\begin{proof}

Recall (Equation \ref{prec}) that for $\gamma \in \cG^y$,
$$
[k_0 \, d_\nu k^{f,l}_1](\gamma) \,\, = \,\, 
\int_{\gamma_1\in \cG^y} k_0(\gamma_1) 
\wh_{\gamma_1}\Bigl{(}d_\nu k^{f,l}_1(\gamma_1^{-1}\gamma)\Bigr{)} dvol_s,
$$
where $\wh_{\gamma_1} :   \nu^*_{s(\gamma_1)}  \to  \nu^*_{r(\gamma_1)}$  is the push-forward map associated to
the holonomy map along the leafwise path $\gamma_1$.  
Set $\what{\gamma}= \gamma_1 \circ \varphi^{f}(f^{-1}s(\gamma_1))$.  Then $k_0(\gamma_1) = k_0^{f,r}(\what{\gamma}) $ (since $f$ is a holonomy diffeomorphism).
As $\what{\gamma}^{-1} = \varphi^{f}(f^{-1}s(\gamma_1))^{-1}\circ \gamma_1^{-1} = 
\varphi^{f^{-1}}(r(\gamma_1^{-1}))\circ \gamma_1^{-1}$,   we also have $d_\nu k^{f,l}_1(\gamma_1^{-1}) \, =  \, ^t(f^{-1}_*)(d_\nu k_1( \what{\gamma}^{-1}))$, and by a trivial extension,
$$
d_\nu k^{f,l}_1(\gamma_1^{-1} \gamma) \, =  \, ^t(f^{-1}_*)(d_\nu k_1( \what{\gamma}^{-1}\gamma)).
$$  
Then   
$$
[k_0 \, d_\nu k^{f,l}_1](\gamma)  \,\, = \,\,  
\int_{\gamma_1\in \cG^y} k_0(\gamma_1) 
\wh_{\gamma_1}\Bigl{(}d_\nu k^{f,l}_1(\gamma_1^{-1}\gamma)\Bigr{)} dvol_s \,\, = \,\,
$$
$$  
\int_{\gamma_1\in \cG^y} k_0^{ f,r}( \what{\gamma})\wh_{\gamma_1}(^t(f^{-1}_*)(d_\nu k_1( \what{\gamma}^{-1}\gamma))) dvol_s \,\, = \,\,  
\int_{\gamma_1\in \cG^y} k_0^{ f,r}( \what{\gamma})\wh_{\gamma_1}(\wh_{\varphi^{f}(f^{-1}s(\gamma_1))}(d_\nu k_1( \what{\gamma}^{-1}\gamma))) dvol_s.
$$ 
As
$\wh_{\gamma_1}\circ \wh_{\varphi^{f}(f^{-1}s(\gamma_1))} = \wh_{\what{\gamma}}  $
and
$\gamma_1 =   \what{\gamma} \circ \varphi^{f^{-1}}(f s(\what{\gamma}))$, this last may be written as
$$
\int_{\what{\gamma} \in \cG^y }  k_0^{ f,r}(\what{\gamma})  \wh_{\what{\gamma}}(d_\nu k_1( \what{\gamma}^{-1}\gamma)) dvol_s   \,\, = \,\, 
[k_0^{f,r} d_\nu k_1](\gamma).
$$
\end{proof}
\end{proof}

As an immediate corollary of Lemma \ref{switchd} we have,
\begin{corollary} \label{extended}
For $k_0,k_1 \in C_c^{\infty}(\cG)$,
$$
k_0 * \delta_\nu k^{f,l}_1 \,\, = \,\,  k_0^{f,r} * \delta_\nu k_1
\quad \text{   and   } \quad
\delta_\nu k_0 * \delta_\nu k^{f,l}_1 \,\, = \,\,   \delta_\nu k_0^{f,r} * \delta_\nu k_1.
$$
\end{corollary}
\begin{proof}
A simple computation gives the first equality.  Applying $\delta_\nu$ and using the facts $\delta_\nu ^2 = 0$ and $\delta_\nu $ is a derivation gives the second.
\end{proof}

Let $g$ be another holonomy diffeomorphism of $(V,F)$.

\begin{proposition}\label{ccoc2} Let $C$ be a closed Haefliger $r$-current and  $\tau_C$ the cyclic cocycle on $C_c^{\infty}(\cG)$ associated with $C$ by Theorem \ref{tauC}.   Then for any
$k_0, \ldots, k_r \in C_c^{\infty}(\cG)$,
$$
\tau_C(k_0, \ldots, k_r \circ \Psi(g)) \,\, = \,\, \tau_C(fk_0, \ldots, fk_{r-1}, (fk_r) \circ \Psi(fgf^{-1})).
$$
\end{proposition}

\begin{proof} 
Note that $\Psi$ is an action, that is $\Psi(gf) = \Psi(g)\Psi(f)$.  If $C$ is a Haefliger $0$-current (so closed as all $0$ currents are) and $k \in C_c^{\infty}(\cG)$, then we have
$$
\tau_C((fk) \circ \Psi(fgf^{-1})) = \tau_C(\Psi(f) \circ k \circ  \Psi(f^{-1})\circ \Psi(fgf^{-1})) =
$$
$$ 
\tau_C(\Psi(f) \circ k \circ \Psi(g) \circ \Psi(f^{-1})) = \tau_C(f ( k \circ \Psi(g))) =
\tau_{f^{-1}_*(C)}( k \circ \Psi(g)) =    \tau_C(k \circ \Psi(g)).
$$
To prove the second to the last equality recall that the action of $f$ preserves the leafwise volume form $dvol_F$. We may assume that the support of $k_1 := k \circ \Psi(g)$ is contained in a fundamental chart for $\cG$, see \cite{BHI}.  As $f$ preserves $\cG^{(0)} = V$, we may in fact assume that $k_1$ is  a smooth compactly supported function on a foliation chart $U \subset V$, with transversal $T_U$.  Then the support of $f(k_1)$ is contained in $W = f(U)$ with transversal $T_W = f(T_U)$. Now 
$\dd \tau_C(f(k_1)) = \int_C \int_F f(k_1) dvol_F$, and $\dd \int_F f(k_1) dvol_F =
\int_W f(k_1)dvol_F$,  so
$$ 
\int_W f(k_1) dvol_F = 
\int_W f(k_1) (f^{-1})^*(dvol_F) = 
\int_W (f^{-1})^*(k_1  dvol_F) = 
\int_{f^{-1}(W)} k_1  dvol_F= 
\int_U k_1 dvol_F.
$$
That is, $\dd \int_F f(k_1) dvol_F  \, |_{T_W} =  \Big(\int_F k_1 dvol_F) \, |_{T_U}\Big) \circ f^{-1}$
and we have the second to the last equality.

To prove the final equality, we need only show that $f^{-1}_*(C) = C$, or equivalently $f_*(C) = C$.  Let $T_U$ be any transversal, and set $T_W = f(T_U)$, also a transversal.  Since $f$ is a holonomy diffeomorphism, we may write $T_U$ as a countable union of open subsets $T_U = \cup U_j$, where $f \, | _{U_j} = h_{\varphi^{f}(x_j)}$ for some $x_j \in U_j$.  Then, $f_*(C) \, |_{f(W_j)} = h_{\varphi^{f}(x_j)*}(C) \, | _{f(W_j)} = C \, | _{f(W_j)}$, and we have this last equality.

Now suppose that $\dim C =r >0$.

\begin{lemma}For $k_1, k_2 \in C_c^{\infty}(\cG)$ and $f \in H$,
$$
fk_1 * \delta_\nu (fk_2) = k_1^{f,l} *  \delta_\nu (k_2^{f^{-1},r}),
$$
and 
$$
 \delta_\nu (fk_1) * \delta_\nu (fk_2) =  \delta_\nu (k_1^{f,l}) *  \delta_\nu (k_2^{f^{-1},r}).
$$
\end{lemma}
\begin{proof}
As above, the second equality is $ \delta_\nu$ applied to the first.

By Corollary \ref{extended}, Lemma \ref{Haction}, and the fact that $\Psi$ is an action, we have 
$$
fk_1 * \delta_\nu (fk_2) =  fk_1 * \delta_\nu ((k_2^{f^{-1},r})^{f,l}) = 
(fk_1)^{f,r} * \delta_\nu (k_2^{f^{-1},r}) = 
$$
$$
((k_1^{f,l})^{f^{-1},r})^{f,r} * \delta_\nu (k_2^{f^{-1},r}) =k_1^{f,l} *  \delta_\nu (k_2^{f^{-1},r}).
$$
\end{proof}

Now
$$
\delta_\nu k^{f^{-1},r}_1  * \delta_\nu k^{f,l}_2 \,\, = \,\, 
\delta_\nu [k^{f^{-1},r}_1  * \delta_\nu k^{f,l}_2] \,\, = \,\,
\delta_\nu [(k^{f^{-1},r}_1)^{f,r}  * \delta_\nu k_2] \,\, = \,\,
\delta_\nu k_1  * \delta_\nu k_2.
$$
By induction it follows immediately that, 
$$
fk_0 * \delta_\nu (fk_1)* \cdots * \delta_\nu (fk_{r-1}) =
k_0^{f,l} * \delta_\nu k_1 * \cdots *   \delta_\nu k_{r-2} *  \delta_\nu (k_{r-1}^{f^{-1},r}).
$$
Similarly, we have 
$$
\delta_\nu (k_{r-1}^{f^{-1},r}) *  \delta_\nu ((fk_{r}) \circ\Psi(fgf^{-1})) =
\delta_\nu k_{r-1} *  \delta_\nu ((k_{r} \circ\Psi(g))^{f^{-1},r}). 
$$
Thus
$$
fk_0 * \delta_\nu (fk_1)* \cdots * \delta_\nu (fk_{r-1}) * \delta_\nu ((fk_{r}) \circ\Psi(fgf^{-1})) =
$$
$$
k_0^{f,l} * \delta_\nu k_1 * \cdots *   \delta_\nu k_{r-1} *  \delta_\nu ((k_{r} \circ\Psi(g))^{f^{-1},r}),
$$
and we have 
$$
\tau_C(fk_0, \ldots, fk_{r-1}, (fk_r) \circ \Psi(fgf^{-1}))  \,\, = \,\, 
\tau_C(k_0^{f,l}, k_1, \ldots, k_{r-1}, (k_{r} \circ\Psi(g))^{f^{-1},r})  \,\, = \,\, 
$$
$$
(-1)^r \tau_C((k_{r} \circ\Psi(g))^{f^{-1},r}, k_0^{f,l}, k_1, \ldots, k_{r-1}),
$$
by the cyclicity of $\tau_C$.
Using  Corollary  \ref{extended}, this equals 
$$
(-1)^r \tau_C(k_{r} \circ\Psi(g), k_0, k_1, \ldots, k_{r-1}) \,\, = \,\,
\tau_C(k_0, \ldots, k_r \circ \Psi(g)),
$$

again by cyclicity.
\end{proof}

We finish this section by specializing to the case of the leafwise action of a compact group $H$. So, we assume that $H$ is a compact group which acts smoothly by diffeomorphisms on the compact manifold $V$.   We assume for simplicity that $H$ preserves the leaves of the foliation $F$.   Then $H$ acts continuously on the algebra $C_c^\infty (\cG, E)$ endowed with its natural compact open topology, as well as on the $C^*$-algebra $C^*(V,F)$,  \cite{Be1}. Given a finite dimensional unitary representation $(X,U)$ of $H$, we denote by ${\underline X}$ the equivariant trivial bundle $V\times X$.

Let $X, X'$ be  finite dimensional Hermitian spaces and $H \to U(X)$, $H\to U(X')$  unitary representations of $H$. As usual, we denote by ${\widetilde{C_c^{\infty}(\cG,E)}}$ the unital algebra $C_c^{\infty}(\cG,E)\oplus \C$. As usual, we denote by ${\widetilde{C_c^{\infty}(\cG,E)}}$ the unital algebra $C_c^{\infty}(\cG,E)\oplus \C$.
Let ${\wtit e} =(e , \Lambda) \in {\widetilde{C_c^{\infty}(\cG,E)}}\otimes \End(X)$ and ${\wtit e}' =(e',\Lambda ') \in {\widetilde{C_c^{\infty}(\cG,E)}}\otimes \End(X')$ be two $H$-invariant  idempotents such that $[\te] - [\te ']$ defines an equivariant $K$-theory class of the algebra $C_c^{\infty}(\cG,E)$.  Recall that $h$ is a holonomy diffeomorphism which generates the compact Lie group $H$. Then we have:

\begin{theorem}\label{locp}  For any closed  Haefliger  current $C$ on $(V, F)$, we have
\begin{itemize}
 \item the complex number 
$$
(\tau_C\sharp \tr)(\Psi^{E\otimes {\underline X}}(h)\circ e, e,...,e) - (\tau_C\sharp \tr)(\Psi^{E\otimes {\underline X}'}(h)\circ e', e',...,e'),
$$
depends only on the equivariant $K$-theory class $[\te] - [\te ']$, and yields an additive map from  the equivariant $K$-theory of $C_c^{\infty}(\cG, E)$ to the scalars.
\item This map induces a pairing
\begin{Equation}\label{pairing2}\
\hspace{3cm}$K^H (C_c^\infty (\cG, E))\otimes H_{ev}(V/F) \longrightarrow C(H)^H.$
\end{Equation}
\item This pairing extends to a well defined pairing
 \begin{Equation}\label{pairing3}\
\hspace{3cm}$K^H (C^*(V,F))\otimes H_{ev}(V/F) \longrightarrow C(H)^H.$
\end{Equation}
\end{itemize}
\end{theorem}

\begin{proof}\
By Proposition \ref{Pairing}, in order to prove the first and second items, we first prove that the functional $C_c^\infty(\cG,E)^{r+1}\times H \to \C$ given by
$$
\phi_C : (k_0, \cdots, k_r | h) \longmapsto \tau_C (\Psi^E(h)\circ k_0, k_1, \cdots, k_r),
$$
is an equivariant continuous cyclic cocycle on $(C_c^\infty(\cG, E), H)$. Continuity is straightforward. We have for any $h'\in H$:
\begin{eqnarray*}
 \phi_C ( h' k_0, \cdots , h' k_r | h' h {h'}^{-1}) & = & \tau_C (\Psi^E (h' h {h'}^{-1}) \circ h'k_0, h'k_1, \cdots, h'k_r)\\
& = & \tau_C (\Psi^E (h')\circ \Psi^E( h)\circ \Psi^E(h')^{-1} \circ h'k_0, h'k_1, \cdots, h'k_r).
\end{eqnarray*}
But recall that $h'k_j = \Psi^E(h') \circ k_j\circ \Psi^E(h')^{-1}$, therefore,
$$
\phi_C ( h' k_0, \cdots , h' k_r | h' h {h'}^{-1})  \,\, = \,\, 
$$
$$
\tau_C (\Psi^E (h')\circ \Psi^E( h)\circ k_0 \Psi^E(h')^{-1} , \Psi^E(h')\circ k_1\circ \Psi^E(h')^{-1}, \cdots, \Psi^E(h')\circ k_r\circ \Psi^E(h')^{-1})  \,\, = \,\, 
$$
$$
\tau_C (\Psi^E (h')\circ \Psi^E( h)\circ k_0  , k_1, \cdots,  k_r\circ \Psi^E(h')^{-1}).
$$
The last equality is a consequence of Proposition \ref{ccoc}. Now by cyclicity of the cochain $\tau_C$ we can again apply Proposition \ref{ccoc} to conclude that
$$
 \phi_C ( h' k_0, \cdots , h' k_r | h' h {h'}^{-1}) = \tau_C (\Psi^E( h)\circ k_0  , k_1, \cdots,  k_r),
$$
and hence the first  and second items.

The third item is a straightforwd consequence of the deep result of Connes \cite{C86}, where he constructs a complicated subalgebra $\maB$ of $C^*(V,F)$ such that
\begin{itemize}
 \item $\maB$ is stable under holomorphic functional calculus in $C^*(V,F)$.
\item $\maB$ contains $C_c^\infty (\cG)$ and so is dense in $C^*(V,F)$.
\item The cyclic cocycle $\phi_C$ on $C_c^\infty (\cG)$  extends to a cyclic cocycle on $\maB$.
\end{itemize}
Hence, the extension of  $\phi_C$ induces an additive map 
$$
\phi_{C, *} : K(\maB) \rightarrow \C,
$$
and the inclusion $i:\maB\hookrightarrow C^*(V,F)$ induces an isomorphism.$K(\maB) \simeq K(C^*(V,F))$.  In the presence of a vector bundle $E$ the same construction yields a subalgebra $\maB^E$ of the twisted by $E$ Connes' $C^*$-algebra $C^*(V,F;E)$, with the same properties.  The main input here is the additional action of the compact Lie group $H$. But, it is obvious from Connes' construction that the algebra $\maB^E$ is an $H$-subalgebra of $C^*(V,F; E)$ and hence, by easy arguments, the $H$-equivariant inclusion $i:\maB^E\hookrightarrow C^*(V,F; E)$ induces an $R(H)$-isomorphism $K^H(\maB^E) \simeq K^H(C^*(V,F;E) (\simeq K^H(C^*(V,F)))$. Moreover, the extension of $\phi_C$ to $\maB^E$ is also clearly an equivariant cyclic cocycle. Therefore, Connes' proof yields the third item in a straightforward way.
\end{proof}

We quote the following proposition for later use. Recall that when the action of $H$ on an algebra $\maA$ is trivial, the equivariant $K$-theory $K^H(\maA)$ is naturally isomorphic to the tensor product $K(\maA)\otimes R(H)$ where $R(H)$ is the representation ring of $H$. 

\begin{proposition}
Assume that the actions of $H$ on $V$ and $E$ are trivial. Then the following diagram is commutative
 
\hspace{0.25in}
\begin{picture}(415,80)

\put(75,60){$ K^H(C_c^\infty (\cG, E))$}
\put(140,50){ $\vector(4,-3){45}$}
\put(125,35){$\left< C, \cdot \right>_h$}

\put(165,70){$\phi$}
\put(150,64){\vector(1,0){35}}

\put(195,60){$K(C_c^\infty (\cG, E)) \otimes R(H)$}
\put(205,50){ $\vector(0,-1){30}$}
\put(203,5){$\C $}
\put(215,35){$\left< C, \cdot \right> \otimes ev_h$}
\end{picture}\\
where $\left< C, \cdot \right>_h$ is the pairing of Theorem \ref{locp} of  the closed Haefliger current $C$ with equivariant $K$-theory, $ev_h: R(H) \to \C$ is evaluation of characters at $h$, and $\phi$ is the natural isomorphism.
\end{proposition}

\begin{proof}
If $\te$ is an idempotent in the unitalization algebra ${\widetilde{C_c^\infty(\cG, E)}}$ and $\rho=(X,U)$ is a finite dimensional unitary representation of $H$, then $\te$ can also be viewed as an equivariant element for the trivial action, and we can therefore consider the equivariant idempotent $\te \rho$ of the algebra ${\widetilde{C_c^\infty(\cG, E)}}\otimes \End(X)$. The map $\te \otimes \rho \mapsto \te \rho$ then implements the isomorphism $\phi^{-1}$. Now, since the action of $H$ on ${\widetilde{C_c^\infty(\cG, E)}}$ is trivial, the multiplier $\Psi^E(h)$ is simply the identity operator. Moreover, the multiplier $\Psi^{E\otimes {\underline X}}(h)$ 
is simply given by $Id_E \otimes U(h)$.   Thus applying the pairing with $C$ to the element $\te \rho$ gives:
\begin{multline*}
 \left<C, [\te \rho]\right>_h = (\tau_C\sharp \tr)(Id \otimes U(h)\circ e\otimes Id_X, e\otimes Id_X,...,e\otimes Id_X)\\ = (\tau_C\sharp \tr)(e\otimes U(h), e\otimes Id_X,...,e\otimes Id_X) = \left<C , [\te]\right> \times \tr (U(h)).
\end{multline*}
Hence the conclusion.
\end{proof}

We now return to the general case, so the action of the compact group $H$ is now a general action by $F$-preserving diffeomorphisms of $(V,F)$.
Recall the prime ideal $I_{[h]}$ associated with (the conjugacy class of) $h\in H$ in $R(H)$:
$$
I_{[h]}:= \{\chi\in R(H) \, | \,  \chi (h)=0\};
$$
The equivariant $K$-theory of the (non $C^*$-) algebra $C_c^\infty(\cG, E)$ is naturally endowed with the structure of an $R(H)$-module.  Therefore we can define the localization of this module at $I_{[h]}$, denoted $K^H(C_c^\infty(\cG,E))_{[h]}$, which is a module  over the localized ring $R(H)_{[h]}$.

\begin{proposition}
Composing  the pairing \ref{pairing2}  with evaluation at the (conjugacy class of the) element $h\in H$, gives the pairing 
$$
\left< \cdot , \cdot \right>_{[h]}: K^H(C_c^\infty(\cG,E)) \otimes H_{ev} (V/F) \longrightarrow \C,
$$
which satisfies 
$$
\left< x \rho , C \right>_{[h]} = \left< x  , C \right>_{[h]} \times \rho(h), \quad x\in K^H(C_c^\infty(\cG,E)), C\in H_{ev} (V/F) \text{ and } \rho\in R(H).
$$
In particular, the pairing $\left< \cdot , \cdot \right>_{[h]}$ induces a pairing of  the Haefliger homology of the foliation with  the equivariant $K$-theory of $C_c^\infty (\cG, E)$ localized at the prime ideal $I_{[h]}$. 
\end{proposition}

\begin{proof}
For simplicity, we will forget the bundle $E$.  Fix two $H$-invariant idempotents $\te = (e, \Lambda) \in C_c^\infty(\cG)\otimes \End(X)$ and $\te'=(e', \Lambda') \in C_c^\infty(\cG)\otimes \End(X')$, and set $x=[\te] - [\te']$. We also consider one (again for simplicity) extra finite dimensional representation $(U'',X'')$ of $H$ which represents a class $\rho$ in $R(H)$. Then the equivariant $K$-theory class $x \rho$ is by definition the class of the formal difference
$$
x\rho = [\te \otimes id_{X''}] - [\te' \otimes id_{X''}].
$$
Therefore, for any even dimensional closed Haefliger  current $C$, we have
\begin{multline*}
\left<x\rho, C \right>_h = (\tau_C\sharp\tr \sharp \tr) ( \Psi^{{\underline X}\otimes {\underline X}''} (h) \circ (e\otimes id_{X''}), e\otimes  id_{X''}, \cdots,  e\otimes  id_{X''}) - \\(\tau_C\sharp\tr \sharp \tr) ( \Psi^{{\underline X}'\otimes {\underline X}''} (h) \circ (e'\otimes id_{X''}), e'\otimes  id_{X''}, \cdots,  e'\otimes  id_{X''}).
\end{multline*}
By definition
$$
\Psi^{{\underline X}\otimes {\underline X}''} (h) \circ (e\otimes id_{X''}) = (\Psi^{{\underline X}}(h)\circ e) \otimes U''(h),
$$
and the same relation holds with $e'$ in place of $e$. Thus 
$$
\left<x\rho, C \right>_h =(\tau_C\sharp\tr) ( \Psi^{{\underline X}} (h) \circ e, e, \cdots,  e) \times \tr (U''(h)) - (\tau_C\sharp\tr) ( \Psi^{{\underline X}'} (h) \circ e', e', \cdots,  e') \times \tr (U''(h)).
$$
As a corollary, we deduce that the formula
$$
\left< \frac{x}{\rho}, C\right>_h := \frac{\left< x, C \right>_h}{ \rho(h)}, \quad x\in K^H(C_c^\infty (\cG)), \,\,  \rho\in R(H)\smallsetminus I_h,
$$
gives the desired pairing between the localized $K$-theory at the conjugacy class of $h$ in $H$, and the Haefliger even dimensional homology of the foliation. The case of general $E$ is similar.
\end{proof}

We now define the higher Lefschetz number. Recall that for any $h\in H$, the topologically cyclic group $H_1$ is the compact (abelian) group topologically generated by $h$ in $H$. Given an $H$-invariant leafwise elliptic complex $(E,d)$ over the foliation $(V, F)$, we can define the localized index class $\frac{\Ind^H (E,d)}{1_{R(H)}}$ with respect to the prime ideal associated with the conjugacy class of $h$ in $R(H)$. So, $\frac{\Ind^H (E,d)}{1_{R(H)}} \in K^H(C^*(V,F))_{[h]}$ and we know that there exists a localized class $\frac{\Ind_\infty^H (E,d)}{1_{R(H)}} \in K^H(C_c^\infty (\cG, E))_h$ which pushes to $\frac{\Ind^H (E,d)}{1_{R(H)}}$ under the functoriality map. Recall that the Lefschetz class of $h$ with respect to $(E,d)$ is given by 
$$
L(h; E, d) = \frac{\Ind^{H_1} (E,d)}{1_{R(H_1)}} \in K^{H_1}(C^*(V,F))_h.
$$
In the same way one defines a smooth Lefschetz class $L_\infty(h;E,d)$ with lives in $K^{H_1}(C_c^\infty(\cG, E))_h$ which pushes to $L(h;E,d)$. One can use naturality with respect to subgroups to check that the following diagram, where $r$ is the natural forgetful map from $H$-equivariant $K$-theory to $H_1$-equivariant $K$-theory, is commutative

\hspace{0.25in}
\begin{picture}(415,80)

\put(75,60){$ K^H(C_c^\infty (\cG, E))$}
\put(140,50){ $\vector(4,-3){45}$}
\put(125,35){$\left< C, \cdot \right>_{[h]}$}

\put(165,70){$r$}
\put(150,64){\vector(1,0){35}}

\put(195,60){$K^{H_1}(C_c^\infty (\cG, E))$}
\put(205,50){ $\vector(0,-1){30}$}
\put(203,5){$\C $}
\put(215,35){$\left< C, \cdot \right>_h$}

\end{picture}\\
Here $C$ is any even dimensional closed Haefliger current on the foliation $(V,F)$ and $\left< C, \cdot \right>_{[h]}$ stands as before for evaluation at $h$ of the pairing of $C$ with the $H$-equivariant $K$-theory, while $\left< C, \cdot \right>_h$ is the same map but for the compact abelian group $H_1$.

Therefore, the Lefschetz numbers given by $h$ evaluted against a given closed Haefliger current $C$ with respect to a given $H$-invariant elliptic complex $(E,d)$ along the leaves of the foliation, do not depend on $H$ as far as this latter exists. We can use any such group and the smallest one is the closed subgroup $H_1$ generated by $h$ in the compact Lie group of isometries of the Riemannian compact manifold $V$. This reduces the Lefschetz problem to the case where $H=H_1$. 

\begin{definition}\ Let $C$ be an even dimensional  closed Haefliger current. Then the $C$-Lefschetz number
of $h$ with respect to a leafwise elliptic pseudodifferential $H$-invariant complex $(E,d)$ is the scalar $L_C(h;E,d)$ defined as
$$
	L_C(h;E,d):=\left< C, L_\infty(h;E,d))\right>_h,
$$
where $L_\infty(h;E,d)$ is the smooth Lefschetz class in the localized equivariant $K$-theory group $K^H(C_c^\infty (\cG, E))_h$.
\end{definition}

As a straightforward consequence of the $K$-theory Lefschetz theorem, we then deduce:
\begin{theorem}\ Let (E,d) be an $H$-invariant leafwise elliptic pseudodifferential complex over
the compact foliated manifold $(V,F)$. Let $C$ be an even dimensional  closed Haefliger current on $(V, F)$. Then we have
$$
L_{C}(h;E,d)\not=0 \quad \Rightarrow \quad  L(h;E,d)\not=0 \quad \Rightarrow \quad \Ind_{V}^{H}(E,d)\not=0 \mbox{ and } V^{h}\not= \emptyset.
$$
\end{theorem}

\begin{proof}\ We know by Proposition \ref{locp} that the map $\left< C, \cdot \right>_{[h]}$  extends to an additive map on the $H$-equivariant $K$-theory $C^*$-algebra $C^*(V,F)$. Since $L(h;E,d) = \frac{\Ind_V^H(E,d)}{1_{R(H)}}$ we immediately deduce, 
$$
L_{C}(h;E,d)\not=0 \quad \Rightarrow \quad L(h;E,d)\not=0 \quad \Rightarrow \quad \Ind_{V}^{H}(E,d)\not=0.
$$
Moreover, by the $K$-theory Lefschetz theorem of \cite{Be1}, we also have
$$
L(h; E,d) \not=0 \quad \Rightarrow \quad V^{h}\not= \emptyset.
$$
\end{proof}

To end this section we explain the case when $(V,F)$ admits a holonomy invariant transverse measure $\Lambda$. Recall that  $\Lambda$ induces a trace on $C^*(V,F; E)$, \cite{CInt}, which is finite on $C_c^{\infty}(\cG,E)$. We will denote this 0-cocycle by $\tau_{\Lambda}$. Because the Ruelle-Sullivan current associated with $\Lambda$ is a Haefliger current we deduce:

\begin{corollary}\label{corlocp}  Let $\te$ be an $H$-invariant projection of
${\widetilde{C_{c}^{\infty}(\cG, E)}} \otimes
\End(X)$, where $(U,X)$ is a finite dimensional unitary representation of $H$. Then the formula:
$$
\tau^{h}_{\Lambda}(\te) = (\tau_{\Lambda} \sharp tr) \left(\Psi^{E\otimes {\underline X}}(h) \circ e \right)
$$
induces an additive map on the (localized at the conjugacy class of $h$) equivariant $K$-theory of the algebra $C_{c}^{\infty}(\cG, E)$:
 $$
\tau^{h}_{\Lambda} : K^{H} (C_{c}^{\infty}(\cG, E))_{[h]}\rightarrow \C.
$$
\end{corollary}

Notice that $\tau^{h}_{\Lambda}$ extends  easily to the localized equivariant $K$-theory of Connes' $C^*$-algebra of the foliation. Corollary \ref{corlocp} was first proved in \cite{Be1}.  Here it is a corollary of Theorem \ref{locp} since $\tau_{\Lambda}$ is an
$H$-equivariant cyclic 0-cocycle. If we define the Lefschetz $\Lambda$-number of $h$ to be
$$
L_{\Lambda}(h ; E,d) =\tau^{h}_{\Lambda}(L_\infty (h ; E,d)),
$$
then we get from the $K$-theory Lefschetz theorem a measured Lefschetz theorem which recovers
the results of \cite{HL1} when the diffeomorphism is an isometry. In this measured case it is easy to see that
$L_{\Lambda}(h;  E,d)$ coincides with the alternating sum of the actions of $h$ on the
kernels of the Laplacians of the
$\cG$-complex so that when $F=T(V)$ we obtain the classical Lefschetz theorem.\\

\section{The higher Lefschetz theorem}\label{HigherLef}

Recall that $(V,F)$ is a closed foliated manifold and that $H$ is a compact Lie group  which 
is generated by a holonomy  diffeomorphism $h$ of $V$, so $h$ is leaf preserving while $H$ is only necessarily $F$ preserving. We assume that $(E,d)$ is a leafwise elliptic pseudodifferential complex on $(V, F)$ which is $H$-equivariant, i.e. $E=\oplus E^i$ and the bundles $E^i$ are $H$-equivariant, while the operators $d^i$ are $H$-invariant. We tackle in this section the main result of this paper, namely the higher Lefschetz theorem in Haefliger cohomology. Rephrasing the results of the previous section and using continuity in $h'\in H$, we have defined an equivariant Chern-Connes character 
$$
\ch^H: K^H (C_c^\infty (\cG, E)) \longrightarrow H_c^{ev} (V/F) \otimes C(H)^H,
$$
which, composed with evaluation at $h$, induces a localized Chern-Connes character 
$$
(\ch^H)_{[h]}: K^H (C_c^\infty (\cG, E))_{[h]} \longrightarrow H_c^{ev} (V/F, \C).
$$
We want to compute the image under $(\ch^H)_{[h]}$ of the Lefschetz class of $h$ with respect to $(E,d)$, in terms of characteristic classes at the fixed points of the elements of $[h]$.  We have adopted the technically simpler choice of applying the pairing of this equivariant Chern-Connes character with closed Haefliger currents. Therefore, we shall rather compute the complex numbers $L_C(h; E,d)$, where $C$ runs through the space of closed Haefliger currents on the foliation $(V, F)$. Due to the naturality of the pairing with respect to subgroups, we can restrict to the case where $H$ is the abelian compact Lie group topologically generated by an element $h$, so the conjugacy class reduces to $h$ and the statements are then much simpler.

We need to restrict Haefliger currents to the fixed point submanifold with its induced foliation, so we  assume that this fixed point submanifold is transverse to the ambiant foliation and has dimension at least equal to the codimension of the foliation, since otherwise the contributions are trivial.  Recall that a smooth submanifold $W$ of the closed manifold $V$ is transverse to the foliation $F$, if the fibers of the tangent space $TW$ to $W$ are transverse as vector spaces to the fibers of $TF$. In the case $\dim (W) \leq \codim (F)$, this means that  $TW$ and $TF$ are in direct sum. In the case $\dim (W) > \codim (F)$ this means that $\dim (TW\cap TF) = \dim (W) - \codim (F)$. We shall then denote by $F_W$ the induced foliation on $W$ whose leaves are the connected components of the intersections of the leaves of $(V, F)$ with $W$. So the bundle $TF_W$ is given by 
$$
TF_W = TW \cap TF.
$$
We shall denote by $\cG_W$ the holonomy groupoid of the foliation $(W, F_W)$. Associated with the inclusion $j: W\hookrightarrow V$ of the transverse submanifold $W$ of $V$, there is a well defined Morita extension map 
$$
j_! : K(C_c^\infty (\cG_W, E|_W)) \longrightarrow K(C_c^\infty (\cG, E))
$$
which is the smooth version of the construction given in \cite{CS}.  When $W$ is in addition an $H$-submanifold, we immediately get an equivariant Morita extension map
$$
j_! : K^H(C_c^\infty (\cG_W, E|_W)) \longrightarrow K^H(C_c^\infty (\cG, E)),
$$
which is defined in exactly the same way. For the convenience of the reader, we recall the construction of the equivariant version of the map $j_!$ below in the proof of Proposition \ref{trans}, see \cite{C82, Be1} for more details.

Let $C$ be any closed Haefliger $k$-current on $(V,F)$ and let $\tau=\tau_C$ be the associated  
$H$-equivariant cyclic $k$-cocycle on the algebra $C_c^{\infty}(\cG, E)$ as in Theorem \ref{tauC}, and recall Theorem \ref{locp}. We assume again that $k$ is even and concentrate on the even $K$-theory. Then we have

\begin{proposition}\label{trans} Let $W$ be a smooth (closed) $H$-submanifold of $V$ which is transverse to the foliation $F$ and has dimension $\geq \codim (F)$.  We  suppose that  $F_W$ is oriented. Then 
\begin{itemize}
 \item There is a well defined $R(H)$-linear morphism $j_! :  K^H(C_{c}^{\infty}(\cG_W, E|_W)) \to K^H(C_{c}^{\infty }(\cG, E))$ which induces the (equivariant) Connes-Skandalis Morita extension morphism \cite{CS}. 
\item The restriction of the current $C$ to the (open) transversals of $F_W$ yields a closed Haefliger current on $(W, F_W)$ that we denote by $C_W$. Moreover, if the generator $h$ preserves the leaves of $F_W$ and is a holonomy diffeomorphism of the foliation $(W, F_W)$, then the following diagram commutes:

\hspace{0.2in}
\begin{picture}(400,80)

\put(85,60){$K^H(C_{c}^{\infty}(\cG_W, E|_W))$}
\put(148,50){ $\vector(4,-3){35}$}
\put(135,30){$\left<C_W, \cdot \right>$}

\put(195,70){$j_!$}
\put(180,64){\vector(1,0){45}}

\put(230,60){$K^H(C_{c}^{\infty }(\cG, E))$}
\put(245,50){ $\vector(-4,-3){35}$}
\put(195,15){$C(H)^H$}
\put(245,30){$\left<C, \cdot \right>$}

\end{picture}
\end{itemize}

\end{proposition}

\begin{proof}\ For the moment, we forget the bundle $E$. Let $\maN$ be an open tubular neighborhood of $W$ in $V$ which is an open $H$-submanifold of $V$ so that its fibres over $W$ are connected open disks of the leaves of $(V,F)$. Notice that since $H$ is compact, such tubular $H$-neighborhoods always exist. As $\maN$ is a transverse submanifold to $(V, F)$, it inherits a foliation $F_\maN$ and a holonomy groupoid $\maG_\maN$ which is clearly an open subgroupoid and submanifold of $\cG$.  As usual,  we identify  $\maN$ with the total space of the normal vector bundle $\pi:N\to W$ to $W$ in $V$. Then the foliation of $\maN$ yields a foliation on $N$ and the action of $H$ on $\maN$ endows $N$ with the structure of an $H$-equivariant vector bundle. More precisely, the vector bundle $N\to W$ is endowed with the smooth foliation whose leaves are the restrictions of $N$ to the leaves of $(W, F_W)$, said differently, it is the foliation on $N$ generated by the integrable subbundle of $TN$ given by the kernel of $\rho\circ \pi_*$ where $\rho: TW \to TW/TF_W$ is the quotient projection. The Morita extension map $j_!$ is then, up to the identification of $(\maN, F_\maN)$ with $(N, F_N)$, the composite map of a Mischenko map induced in $K^H$-theory by an algebra morphism
$$
M: C_{c}^{\infty}(\cG_W) \longrightarrow C_{c}^{\infty}(\cG_N),
$$
that we recall below, and the excision map induced by the trivial extension $\iota: C_{c}^{\infty}(\cG_\maN)\to C_{c}^{\infty}(\cG)$ associated with the open subgroupoid $\cG_\maN$ of $\cG$. 
Recall that $\nu$ is the transverse bundle to the ambiant foliation.  We also denote by $\nu$  its restriction to $\maN$ and $W$, as well as  the normal bundle to the foliation of the vector bundle $N$.   It should be clear from the context which we mean.

We first show that the equivariant pairing commutes with the equivariant excision.  Any smooth compactly supported form $\omega\in C_c^\infty (\cG_\maN, \wedge^jr^*(\nu^*))$ can be trivially extended to a smooth compactly supported form $\iota_*\omega \in  C_c^\infty (\cG, \wedge^jr^*(\nu^*))$. Moreover, we get by direct inspection that for any $\omega_1, \omega_2\in C_c^\infty (\cG_\maN, \wedge^jr^*(\nu^*))$,
$$
(\iota_*\omega_1) (\iota_*\omega_2) = \iota_*(\omega_1\omega_2) \quad \text{ and } \quad 
d_\nu \circ \iota_* = \iota_*\circ d_\nu.
$$
On the other hand, notice that since $h$ preserves the leaves of $F_W$ and is a holonomy diffeomorphism of $(W, F_W)$, it also preserves the leaves of $F_\maN$ and is a holonomy diffeomorphism of $(\maN, F_\maN)$. Hence, for any $h\in H$, the multiplier $\Psi (h)$ can be defined as an operator on $L^2((\cG_\maN)_x)$.   Then 
$$
(\Psi(h)\circ \iota k_0) d_\nu (\iota k_1) \cdots d_\nu (\iota k_r) = \iota_* \left[(\Psi(h)\circ  k_0) (d_\nu  k_1) \cdots (d_\nu k_r)\right],
$$
and a simple algebraic computation shows that 
$$
\left[(\Psi(h)\circ \iota k_0)* \delta_\nu (\iota k_1)* \cdots *\delta_\nu (\iota k_r) \right]_{11}= \iota_* \left[(\Psi(h)\circ  k_0)* (\delta_\nu  k_1)* \cdots * (\delta_\nu k_r)\right]_{11}.
$$
Integration over the leaves also commutes with $\iota_*$ and we deduce 
$$
\tau_C (\Psi(h)\circ \iota k_0, \iota k_1, \cdots, \iota k_r) = \tau_{C_\maN} (\Psi(h) \circ k_0, k_1, \cdots, k_r), \quad k_j\in C_c^\infty (\cG_\maN).
$$
This finishes the proof that the equivariant pairing with $C$ commutes with the excision map. 

The equivariant leafwise diffeomorphism between $(\maN, F_\maN)$ and $(N, F_N)$ induces an isomorphism in equivariant $K$-theory, so we now concentrate on the Mishenko map 
$M:C^{\infty}_c(\cG_W) \to C^{\infty}_c(\cG_N)$ associated with the vector bundle $\pi: N\to W$. We shall identify $W$ with the zero section of the vector bundle $(N, \pi, W)$. It is obvious that  the holonomy groupoid $\maG_N$ may be identified with the smooth groupoid whose space of arrows is 
$$
\{(\gamma, n, n')\in \cG_W\times N^2 \text{ such that } \pi(n)=s(\gamma)\text{ and } \pi (n')=r(\gamma)\}.
$$
The groupoid laws should be clear.   Choose a Lebesgue class measure $\beta$ along the fibers of $(N,  \pi, W)$ which is given by a volume form along the fibers, and a Lebesgue class measure $\mu$ along the leaves of $F_W$ associated with a leafwise volume form. Then, there is a Lebesgue class measure along the leaves of $F_N$ which is locally given by $\mu\otimes \beta$ and that is associated with a volume form along the leaves of $F_N$. Using a partition of unity argument and averaging if necessary,  construct a non negative smooth compactly supported function $\varphi$ on the total space $N$ which is $H$-invariant and such that for every $w\in W$,  the restriction of $\varphi$ to the fiber $N_w$ has $L^2$-norm equal to $1$, i.e. 
$$
\int_{N_w} \varphi(n)^2 d\beta (n) = 1, \quad \forall w\in W.
$$
The homomorphism $M$ applied to $k \in C^{\infty}_c(\cG_W)$ is then simply given by:
$$
M(k) (\gamma, n, n') := k(\gamma) \varphi (n) \varphi(n').
$$
It is easy to check that $M$ is an algebra morphism which is $H$-equivariant. With obvious notations, the bundle $\pi_*\nu$ is a normal bundle to the foliation $F_W$ and coincides with the restriction of $\nu$ to $W$.   
Let $[\alpha]= [(\gamma, n, n')]\in \cG_N$, so $\alpha$ is a path joining $n$ to $n'$ inside a leaf of $F_N$.  Denote by ${h}_\alpha$ the holonomy local diffeomorphism of $\alpha$ acting between a small transversal at $n$ and a small transversal at $n'$.  
Denote by ${h}_\gamma$ the local holonomy diffeomorphism associated with the leafwise (in $F_W$) path $\gamma = \pi(\alpha)$. Notice that, up to the identifications through $\pi$ of the transversals at points of $N$ and their projected transversals in $W$, the local diffeomorphism ${h}_\alpha$ coincides with the local diffeomorphism ${ h}_\gamma$,  that is
$$
\pi \circ { h}_\alpha = {h}_\gamma \circ {\pi}.
$$
Therefore at the level of the induced tangent maps, we deduce a similar relation in the actions on the normal bundle $\nu$. This shows, by definition of the transverse differential $d_\nu$ corresponding to the choice of $\nu$, that for any $X\in \nu_{n'}$, and denoting by $Y\in T_{[\alpha]}\cG_N$ the unique lift with $r_*Y = X$ and $s_*Y= {h}^{-1}_{\alpha, *} X$, we have 
$$
< [d_\nu (M (k))]_{[\alpha]} , X>  \,\, = \,\, < d (M(k))_{[\alpha]}, Y >   \,\, = \,\,  
$$
$$\varphi(n) \varphi(n') <d_\nu k, Y>  \,\, + \,\,  \varphi(n) k(\gamma) <d_{n'}\varphi, X>  \,\, + \,\,  \varphi(n') k(\gamma) <d_n\varphi, {h}^{-1}_{\alpha,*} (X)>.
$$
Since $\varphi$ is smooth compactly supported, we have the following equality in $\nu^*_w$ for any $w\in W$,
$$
\int_{N_w}  (\pi^*_n)^{-1} [\varphi(n) (d_{\nu} \varphi)(n)] d\beta (n) \,\, = \,\, \frac{1}{2} d_{\nu} \left[\int_{N_w} (\varphi(n))^2 d\beta(n) \right] = 0,
$$
where $\pi^*_n:\nu^*_w  \simeq \nu^*_n$.
Now consider the contribution of the second term in $d_\nu (M (k_1))$ to say $M(k_0) d_{\nu}(M(k_1))$, evaluated at the unit $[n']$ of $\cG_N$ determined by the point $n' \in N$.  Recall that for $ \alpha = (\gamma, n, n')$, $\widetilde{h}_{\alpha} =  {^t{ h}_{\alpha^{-1}, *}}:\nu^*_{n} \to \nu^*_{n'}$, and $\widetilde{h}_{\gamma} =  {^t{ h}_{\gamma^{-1}, *}}:\nu^*_{\pi(n)} \to \nu^*_{\pi(n')}$.  Using the above equality we have,
\begin{eqnarray*}
 \left(M(k_0) [\varphi \pi^*(k_1) d_\nu \varphi]\right) ([n'])  & = &   \int_{[\alpha]\in \cG_N^{n'}} M(k_0) (\alpha)\widetilde{h}_{\alpha} ([\varphi \pi^*(k_1) d_\nu \varphi](\alpha^{-1})) d(\mu\otimes \beta) ([\alpha])  \\ & = & \int_{\gamma\in \cG_W^{\pi(n')}} \int_{n\in N_{s(\gamma)}}  k_0(\gamma) \varphi(n') \varphi(n) k_1(\gamma^{-1}) \varphi(n') \widetilde{h}_{\alpha} ( d_\nu \varphi (n)) d\beta(n) d\mu(\gamma)\\ & = &  \varphi(n')^2 \int_{\gamma \in \cG_W^{\pi(n')}} k_0(\gamma) k_1(\gamma^{-1}) \left[\int_{n\in N_{s(\gamma)}} \widetilde{h}_{\alpha} (\varphi(n) d_\nu \varphi (n)) d\beta(n)\right] d\mu(\gamma).
\end{eqnarray*}
But 
\begin{eqnarray*}
 \int_{n\in N_{s(\gamma)}}\widetilde{h}_{\alpha} (\varphi(n) d_\nu \varphi (n)) d\beta(n) & = &
\pi^*_{n'} \int_{n\in N_{s(\gamma)}} (\pi^*_{n'})^{-1} \widetilde{h}_{\alpha} (\varphi(n) d_\nu \varphi (n)) d\beta(n)   
\\ & = & \pi^*_{n'}  \int_{n\in N_{s(\gamma)}} \widetilde{h}_{\gamma}(\pi^*_{n})^{-1} (\varphi(n) d_\nu \varphi (n)) d\beta(n)\\ & = & \pi^*_{n'}\widetilde{h}_{\gamma} \int_{n\in N_{s(\gamma)}} (\pi^*_{n})^{-1} (\varphi(n) d_\nu \varphi (n)) d\beta(n) \,\, = \,\,  0.
\end{eqnarray*}
Essentially the same computation shows that, restricted to the units,  $M(k_0)\times $the third term in $d_\nu (M (k_1))$ is also zero. Thus $[M(k_0) d_\nu (M(k_1))] ([n'])= \varphi(n')^2  [k_0 d_\nu k_1] ([\pi(n')])$.   A similar (more involved but straightforward) computation gives that   restricted to the units, all the terms in the differential form
$$
M(k_0) d_{\nu} (M(k_1)) \cdots d_\nu (M(k_r))
$$
which involve at least one differential of $\varphi$ are zero.  Hence 
$$
\left[M(k_0) d_{\nu} (M(k_1)) \cdots d_\nu (M(k_r)) \right]([n']) = 
\varphi(n')^2 \left[k_0 d_\nu(k_1) \cdots d_\nu(k_r)\right]([\pi(n')]).
$$
In order to extend this relation to the modified differential graded algebras, we note that the argument above would suffice if the curvature $\theta$ evaluated at $[(\gamma, n, n')]$  depended only on $\gamma$, so that it would pull out of the integration just as $k_j(\gamma)$ does.   We show that in fact this is morally true in the sense that we may treat $\theta$ as if this were true with no ill effects.  First note that the Haefliger class we pair with is independent of the choice of the function $\varphi$.  Now consider the one parameter family of such functions given by $\varphi^2_t(n) = t^{-k}\varphi^2(n/t)$. Here $k$ is the fiber dimension of $N$.
Next note that for the general case, we must deal with integrals of the form 
$$
\int_{n\in N_{s(\gamma)}}  \varphi_t(n) d_\nu \varphi_t (n) \theta(\gamma, n, n') d\beta(n)  \quad \text{and} \quad \int_{n\in N_{s(\gamma)}} \varphi^2_t(n) \theta(\gamma, n, n')  d\beta(n).
$$
For simplicity, we are ignoring things such as $\widetilde{h}_{\alpha}$,  $ \pi^*_{n'}$, $(\pi^*_{n})^{-1}$, etc.  
Now we may write $\theta(\gamma, n, n') =  \theta(\gamma, 0, n') + \sum_j n_j \theta_j(\gamma, n, n')$, where the $\theta_j$ and all their derivatives are uniformly bounded on $\cG$.   Then, by a simple change of variables, we have
$$
\int_{n\in N_{s(\gamma)}} \varphi^2_t(n) \theta(\gamma, n, n')  d\beta(n) \,\, = \,\,  \theta(\gamma, 0, n') + \sum_j     \int_{n\in N_{s(\gamma)}} \varphi^2(n) t n_j \theta_j(\gamma, tn, n')   d\beta(n)
\,\, = \,\,  \theta(\gamma, 0, n') + \Phi_t.
$$
Since the $\theta_j$ and all their derivatives are uniformly bounded on $\cG$, $\Phi_t$ is smooth in all its variables, and all its transverse derivatives converge to zero as $t$ goes to zero.  Since $d_{\nu} \theta = 0$, twice the first integral is just $d_\nu$ applied to the second integral, so equals $d_\nu \Phi_t$, and behaves just like $\Phi_t$.   Thus we get
$$
\left[M(k_0) * \delta_{\nu} (M(k_1)) * \cdots * \delta_\nu (M(k_r)) \right]_{11}([n']) = 
\varphi(n')^2 \left[k_0 * \delta_\nu(k_1) * \cdots * \delta_\nu(k_r)\right]_{11}([\pi(n')]) + \Psi_t,
$$
where $\Psi_t$ is smooth in all its variables, and all its transverse derivatives converge to zero as $t$ goes to zero.  By Fubini,  integration over the leaves of $F_N$ becomes integration over the fibers of $N$ followed by integration over the leaves of $F_W$.   Thus, up to the identification of the Haefliger forms on $(W, F_W)$ with the Haefliger forms on $(N, F_N)$, 
$$
\int_{F_N} \left[M(k_0) * \delta_{\nu} (M(k_1)) * \cdots * \delta_\nu (M(k_r)) \right]_{11} = \int_{F_W} \left[k_0 * \delta_\nu(k_1) * \cdots * \delta_\nu(k_r)\right]_{11} +  \int_{F_N}  \Psi_t.
$$
Now, integration over the leaves of $F_N$ is really integration over compact sets, and since the integrands are uniformly bounded,  this integration commutes with taking transverse derivatives. The Haefliger class determined by $\dd \int_{F_N}  \Psi_t$ is independent of $t$, and letting $t \to 0$, we see that it has a representative  whose derivatives up to any finite order are as small as we please.  Thus it is
the zero Haefliger class, and so it contributes nothing to the pairing.

To finish the proof of the commuting with the equivariant Mischenko map, notice that since the projection $\pi$ is $H$-equivariant, we have for any $h\in H$ and any $n\in N$, ${\widetilde\pi} (\varphi^h (n)) = \varphi^h (\pi(n))$, where ${\widetilde\pi}: \cG_N\to \cG_W$ is the projection.  Hence
$$
\Psi(h)\circ M(k) = M(\Psi(h)\circ k), \quad \forall k\in C_c^\infty(\cG_W).
$$

 We now indicate how to handle the inclusion of the bundle $E$, and for this we revert to  the 
open tubular neighborhood $\maN$.  Denote by $\rho:\maN \to W$ the projection. 
In order to extend the above construction, we only need a smooth way to identify $E_x$ with $E_{\rho(x)}$ for each $x \in \maN$.   Choose a connection on $E$.  We may assume that the fibers of $\maN$ as so small that there is a unique (leafwise for $F$) geodesic from $x$ to $\rho(x)$ in fiber of $x$.  We then use the parallel translation defined by the connection on $E$ along this geodesic to identify $E_x$ with $E_{\rho(x)}$.
\end{proof}

We are now in position to state the following

\begin{theorem}\label{basic}  [Higher Lefschetz Theorem]
Let $F$ be an oriented foliation of the Riemannian manifold $(V,g)$.  Let $h$ be a holonomy isometry of $(V,g)$ (so $h$ preserves the leaves of $F$).  Denote by $H$ the compact Lie group generated by $h$ in the group of isometries of $(V,g)$.  Assume that the fixed-point submanifold $V^h=V^H$ of $h$ is transverse to the foliation and denote by $F^h$ its induced foliation and by $N^h$ the normal bundle to $V^h$ in $V$. Then for any leafwise elliptic $H$-equivariant pseudodifferential complex $(E,d)$  over 
$(V,F)$ and
every closed even dimensional Haefliger current $C$ on $(V,F)$, the higher Lefschetz number of $h$ with respect to $(E,d)$ is given by
$$
	L_{C}(h;E,d) \,\, = \,\, \Ind_{C|_{V^h}} \left(\frac{i^{*}[\sigma(E,d)](h)}
	{\lambda_{-1}(N^{h}\otimes{\C})(h)} \right),
$$
where $C|_{V^h}$ is the closed Haefliger current on $(V^h, F^h)$ which is the restriction of $C$, $\Ind_{C|_{V^h}}:K(TF^{h})\otimes \C \rightarrow \C$ is  the complexified higher $C|_{V^h}$-index map on $(V^h, F^h)$, see \cite{BHI}, and $i:TF^{h} \hookrightarrow TF$ is the H-inclusion.
\end{theorem}

\begin{proof}\  
We will use Theorem \ref{abgrp}, which can be summarized by the
commutativity of the following square\\

\hspace{0.2in}
\begin{picture}(415,90)

\put(105,75){$ K_{H}(TF^h)_h$}
\put(125,65){ $\vector(0,-1){35}$}
\put(80,48){$(\Ind_{V^{h}}^H)_h$}

\put(195,85){${(i_!)}_h$}
\put(160,79){\vector(1,0){65}}

\put(195,25){${(j_!)}_h$}
\put(180,19){\vector(1,0){45}}

\put(235,75){ $K_{H}(TF)_h$}
\put(255,65){ $\vector(0,-1){40}$}
\put(270,48){$(\Ind^{H}_{V})_h$}

\put(75,15){$K^H(C^*(V^h, F^h))_h$}

\put(235,15){$K^H(C^*(V,F))_h$}

\end{picture} \\
and the fact that  $(i_!)_{h}$ is an $R(H)_{h}$-isomorphism with inverse given by $\dd \frac{i_{h}^{*}}{\lambda_{-1}(N^{h}\otimes{\C})}$.  
Recall that the Morita extension $j_!$ is induced by the smooth Morita extension $j_!:K^H(C_c^{\infty}(\cG_{V^h}, E_{V^h})) \to K^H(C_c^{\infty}(\cG, E))$ described in the previous proposition applied to  the transverse submanifold $W=V^h$. We deduce that the Lefschetz class $L(h ; E,d):= \frac{\Ind_V^H [\sigma(E,d)]}{1_{R(H)}}$ is given by
$$
	L(h ; E,d) \,\, = \,\,  (j_!)_{h} \circ (\Ind_{V^h} \otimes R(H)_h)
	\left(\frac{i^{\ast}[\sigma(E,d)]}{\lambda_{-1}
	(N^{h} \otimes \C)}\right).
$$
On the other hand, $L_{C}(h;E,d)$ is the image of $L(h;E,d)$ by the additive map $<C, \cdot>_h$, so
$$
	L_{C}(h;E,d)  \,\, = \,\, (<C, \cdot>_h \circ j_!)_{h} \circ (\Ind_{V^h}^H)_{h})
	\left(\frac{i^{\ast}[\sigma(E,d)]}{\lambda_{-1}(N^{h} \otimes \C)}\right).
$$
By Proposition \ref{trans}, we know that 
$$
	<C, \cdot>_h \circ (j_!)_h  \,\, = \,\, <C|_{V^h}, \cdot >_h.
$$
Since the action of $H$ on $V^h$ is trivial, the $H$-algebra $C^*(V^h, F^h)$ is also endowed with the trivial action. Hence there is a natural $R(H)_h$-isomorphism 
$$
K^H(C^*(V^h, F^h))_h \,\,  \simeq \,\,  K(C^*(V^h, F^h)) \otimes R(H)_h,
$$
with respect to which the additive map $<C|_{V^h}, \cdot >_h$ reads simply  $<C|_{V^h}, \cdot > \otimes ev_h$, where 
$$
<C|_{V^h}, \cdot > : K(C_c^\infty (C^*(V^h, F^h)) \longrightarrow \C,
$$
is the non equivariant pairing with the closed Haefliger current $C|_{V^h}$ and where $ev_{h}:R(H)_{h} \rightarrow \C$ is evaluation at $h$, given by $ev_{h}(\chi/\rho)=\chi(h)/\rho(h)$.
So, we get
$$
L_{C}(h;E,d) \,\, = \,\, ([<C|_{V^{h}}, \cdot > \circ \Ind_{V^h}]\otimes ev_{h})
	\left(\frac{i^{\ast}[\sigma(E,d)]}{\lambda_{-1}(N^{h} \otimes \C)}\right) = (\Ind_{C|_{V^h}} \otimes ev_{h})\left(\frac{i^{\ast}[\sigma(E,d)]}{\lambda_{-1}(N^{h} \otimes \C)}\right).
$$ 
where $\Ind_{C|_{V^h}}= \, <C|_{V^{h}}, \cdot > \circ \Ind_{V^h}$ is the complexified index map associated with $C|_{V^h}$ on the foliation of the fixed point submanifold $V^h$.

In the same way, we  define evaluation at $h$
$$
	K_{H}(TF^{h})_{h}  \rightarrow K(TF^{h}) \otimes \C,
$$
which we denote by $x \mapsto x(h)$. More precisely, we use again that $K_H(TF^h)_h \simeq K(TF^{h}) \otimes R(H)_{h}$ and evaluate at $h$ the elements of $R(H)_h$. We finally get
$$
	L_{C}(h;E,d) \,\, = \,\,\Ind_{C|_{V^h}} \left(\frac{i^{*}[\sigma(E,d)](h)}{\lambda_{-1}(N^h\otimes{\C})(h)} \right)
$$
as claimed.
\end{proof}

We finish this section by stating the cohomological Lefschetz formula obtained using the higher index theorem for Haefliger currents \cite{BHI, C94}.

\begin{theorem}\label{basic2}\ [Cohomological Lefschetz formula]
Under the assumptions of Theorem \ref{basic} and when the foliation $F^h$ of the fixed point submanifold $V^h$ is oriented,   we get for any even dimensional closed Haefliger current $C$, 
$$
	L_C(h; E,d)=\left<\int_{F^h} \frac{ch_{\C}(i^{*}[\sigma(E,d)](h))}{ch_{\C}(\lambda_{-1}(N^{h}
	\otimes{\C})(h))} \, \Td(F^{h} \otimes \C) , C|_{V^h} \right >
$$
where $\Td$ denotes the Todd characteristic class of complex bundles and $\int_{F^h}$ is the integration over the leaves of the foliation $(V^h, F^h)$.
\end{theorem}

\begin{proof}\ Let $\phi_{C, V^{h}}$ be the map from $H_c^*(TF^h, \R)$  to $\C$ given by
$$ 
	\phi_{C, V^{h}}(x) = \left<\int_{F^h} [x \, \Td(TF^{h} \otimes \C)] , C|_{V^{h}} \right >.
$$
Then the higher index theorem for foliations, applied to the closed foliated leafwise oriented manifold $(V^h, F^h)$ with the closed even dimensional Haefliger current  $C|_{V^{h}}$, reduces (up to a sign that we include in the definition of the Chern character), to the  equality
$$
	\Ind_{C, V^{h}} = \phi_{C, V^{h}} \circ \ch,
$$
where $\ch:  K(TF^h) \to H_c^*(TF^h, \R)$ is the usual Chern character for compactly supported $K$-theory.
This Chern character can be extended trivially to $K(TF^{h}) \otimes R(H)_{h}$ and then composed with evaluation at $h$ to yield 
 the map $\theta: K(TF^{h}) \otimes R(H)_{h} \to H^{*}(TF^{h},\C)$ which is given by
$$
	\theta(x\otimes {\frac{\chi}{\rho}}) \,\,= \,\,  \ch_{\C} \left [ (x \otimes \chi / \rho) (h) \right ]\,\,= \,\, 
	\ch_{\C} \left [ (x \otimes \chi(h) / \rho(h))  \right ] \,\,= \,\, [\chi(h) / \rho(h)] \, \ch(x).
$$
Hence, if we trivially extend $\phi_{C, V^{h}}$ to $H^{*}(TF^{h},\C)$, we get
\begin{eqnarray*}
	(\phi_{C, V^{h}} \circ \theta) \left ( x \otimes \chi / \rho \right )& = &
	\left < \ch(x) \, \frac {\chi(h)}{\rho(h)} \,  \Td(TF^{h} \otimes \C) , C|_{V^{h}} \right >\\
	 &=&  \left < \ch_{\C}((x \otimes \frac{\chi}{\rho})(h)) \, \Td(TF^{h} \otimes \C) , C|_{V^{h}} \right >.
\end{eqnarray*}
On the other hand, by the  $K$-theory Lefschetz theorem, we have
$$
	L_C(h;E,d) = \Ind_{C, V^{h}} \left ( \frac{i^{*}[\sigma(E,d)](h)}{\lambda_{-1}(N^{h}\otimes{\C})(h)} \right ).
$$
As
$$ 
	(\phi_{C, V^{h}} \circ \theta)(u)=\Ind_{C, V^{h}}(u(h)),
$$
we obtain
$$
	L_C(h;E,d) = (\phi_{C, V^h} \circ \theta) \left ( \frac{i^{*}[\sigma(E,d)](h)}
	{\lambda_{-1}(N^{h}\otimes{\C})(h)} \right ),
$$
and finally
$$
	L_C(h;E,d)=\left< \frac{ch_{\C}(i^{*}[\sigma(E,d)](h))}
	{ch_{\C}(\lambda_{-1}(N^{h}\otimes{\C})(h))} \, \Td(TF^{h} \otimes \C) , C|_{V^{h}} \right >
$$
as claimed.
\end{proof}

Theorem \ref{basic2} simplifies notably when the fixed point 
submanifold $V^{h}$ is a strict transversal, that  is when $V^h$ is transverse to the foliation
with dimension exactly the codimension of the foliation. This simplification corresponds in the case of a foliation with one leaf  to 
the original case of isolated fixed points.  In particular,
\begin{corollary}\  \label{sttrans}
Suppose that $V^{h}$ is a strict transversal and $C$ is an even  closed  Haefliger current.   Then under the asumptions of Theorem \ref{basic}, 
$$
L_C(h;E,d)  \,\, = \,\,  \left < \frac{ \sum_{i} (-1)^{i} \ch_{\C}([E^{i}|_{V^{h}}](h))}{\sum_{j} (-1)^{j} \ch_{\C}
([\wedge^{j}(TF|_{V^{h}} \otimes \C)](h))} , C|_{V^h} \right >.
$$
\end{corollary}

\begin{proof}\  Here $TF^{h} \cong V^{h}$, so
$$
i^{*}[\sigma(E,d)]=\sum_{i} (-1)^{i} [E^{i}|_{V^{h}}] \mbox{, } N^{h}\cong TF|_{V^{h}} \mbox{ and } 
\Td(TF^{h} \otimes \C)=\Td(V^{h}\times \C)=1.
$$
Applying Theorem \ref{basic2} gives the corollary.. 
\end{proof}

We point out that the term $\ch_{\C}(i^{*}[\sigma(E,d)](h))$ is 
not easy to compute in  general.   The simplifications that occur in the most important geometric 
cases, together with some applications of the higher Lefschetz formula, are treated in the second part of this work \cite{BHV}.

\end{document}